\documentclass[12pt,psamsfonts,leqno,oneside,letterpaper]{amsart}
\usepackage[dvips,text={6.5truein,9truein},left=1truein,top=1truein]{geometry}
\usepackage{amssymb,amsmath,amscd,enumerate,
}
\usepackage[pdftex]{graphicx}
\usepackage{url}

\usepackage[colorlinks,linkcolor=blue,citecolor=blue,pdfstartview=FitH]{hyperref}
\input xy
\xyoption{all}
\SelectTips{cm}{12}

\usepackage{color}
\usepackage{mathrsfs}
\newcommand\invisiblecomment[1]{\empty}

\newcommand\missingref[1]{\empty}
\newcommand\tinymissingref[1]{\empty}
\newcommand\abstractcomment[1]{\empty}

\theoremstyle{definition}
\newtheorem{para}{}[section]

\newtheorem{remark}[para]{Remark}
\newtheorem{reformulation}[para]{Reformulation}
\newtheorem{remarks}[para]{Remarks}
\newtheorem{notation}[para]{Notation}

\newtheorem{convention}[para]{Convention}
\newtheorem{definition}[para]{Definition}
\newtheorem{definitions}[para]{Definitions}
\newtheorem{definitionnotation}[para]{Definition and Notation}
\newtheorem{notationdefinitionremark}[para]{Notation, Definition and Remark}

\newtheorem{remarksnotation}[para]{Remarks and Notation}

\newtheorem{remarknotation}[para]{Remark and Notation}
\newtheorem{notationremark}[para]{Notation and Remark}
\newtheorem{notationreviewremarks}[para]{Notation, Review and Remarks}
\newtheorem{definitionremark}[para]{Definition and Remark}
\newtheorem{definitionsremarks}[para]{Definitions and Remarks}
\newtheorem{notationremarks}[para]{Notation and Remarks}
\newtheorem{definitionsnotation}[para]{Definitions and Notation}

\newtheorem{reviewdefinition}[para]{Review and Definition}

\newtheorem{definitionnotationremarks}[para]{Definition, Notation and Remarks}
\newtheorem{definitionsnotationremarks}[para]{Definitions, Notation and Remarks}

\newcommand\Alternatives{\begin{enumerate}[(i)]}
\newcommand\EndAlternatives{\end{enumerate}}
\newcommand\Conditions{\begin{enumerate}[(1)]}
\newcommand\EndConditions{\end{enumerate}}

\theoremstyle{plain}
\newtheorem{theorem}[para]{Theorem}

\newtheorem{lemma}[para]{Lemma}
\newtheorem{remarkdefinition}[para]{Remark and Definition}
\newtheorem{remarkdefinitions}[para]{Remark and Definitions}
\newtheorem{proposition}[para]{Proposition}

\newtheorem{corollary}[para]{Corollary}
\newtheorem{conjecture}[para]{Conjecture}

\newtheorem*{introclosedhomologytheorem}{Theorem 7.8}
\newtheorem*{intromaybecuspstheorem}{Theorem 7.9}
\newtheorem*{intromainertheorem}{Theorem 7.6}
\newtheorem*{introsemimainertheorem}{Theorem 7.7}

\newtheorem{claim}[equation]{}
\numberwithin{equation}{para}
\numberwithin{figure}{section}
\numberwithin{specialremark}{para}

\newcommand\Number{\begin{para}}
\newcommand\EndNumber{\end{para}}
\newcommand\Definition{\begin{definition}}
\newcommand\EndDefinition{\end{definition}}
\newcommand\Definitions{\begin{definitions}}
\newcommand\DefinitionsNotation{\begin{definitionsnotation}}
\newcommand\NotationDefinitionRemark{\begin{notationdefinitionremark}}
\newcommand\EndNotationDefinitionRemark{\end{notationdefinitionremark}}

\newcommand\DefinitionNotation{\begin{definitionnotation}}
\newcommand\RemarksNotation{\begin{remarksnotation}}
\newcommand\Reformulation{\begin{reformulation}}
\newcommand\EndRemarksNotation{\end{remarksnotation}}
\newcommand\EndReformulation{\end{reformulation}}
\newcommand\RemarkNotation{\begin{remarknotation}}
\newcommand\EndRemarkNotation{\end{remarknotation}}
\newcommand\NotationRemark{\begin{notationremark}}
\newcommand\EndDefinitionNotationRemarks{\end{definitionnotationremarks}}
\newcommand\NotationReviewRemarks{\begin{notationreviewremarks}}
\newcommand\DefinitionRemark{\begin{definitionremark}}
\newcommand\DefinitionsRemarks{\begin{definitionsremarks}}
\newcommand\DefinitionNotationRemarks{\begin{definitionnotationremarks}}
\newcommand\DefinitionsNotationRemarks{\begin{definitionsnotationremarks}}
\newcommand\EndDefinitionsNotationRemarks{\end{definitionsnotationremarks}}

\newcommand\NotationRemarks{\begin{notationremarks}}
\newcommand\EndNotationRemark{\end{notationremark}}
\newcommand\EndNotationReviewRemarks{\end{notationreviewremarks}}
\newcommand\EndDefinitionRemark{\end{definitionremark}}
\newcommand\EndDefinitionsRemarks{\end{definitionsremarks}}
\newcommand\EndNotationRemarks{\end{notationremarks}}
\newcommand\EndRemarkDefinition{\end{remarkdefinition}}
\newcommand\EndRemarkDefinitions{\end{remarkdefinitions}}
\newcommand\RemarkDefinition{\begin{remarkdefinition}}
\newcommand\RemarkDefinitions{\begin{remarkdefinitions}}
\newcommand\EndDefinitionsNotation{\end{definitionsnotation}}
\newcommand\EndDefinitionNotation{\end{definitionnotation}}
\newcommand\ReviewDefinition{\begin{reviewdefinition}}
\newcommand\EndReviewDefinition{\end{reviewdefinition}}
\newcommand\EndDefinitions{\end{definitions}}
\newcommand\Theorem{\begin{theorem}}
\newcommand\EndTheorem{\end{theorem}}
\newcommand\Conjecture{\begin{conjecture}}
\newcommand\EndConjecture{\end{conjecture}}
\newcommand\Remark{\begin{remark}}
\newcommand\EndRemark{\end{remark}}
\newcommand\Remarks{\begin{remarks}}
\newcommand\EndRemarks{\end{remarks}}
\newcommand\Convention{\begin{convention}}
\newcommand\EndConvention{\end{convention}}
\newcommand\Notation{\begin{notation}}
\newcommand\EndNotation{\end{notation}}
\newcommand\Lemma{\begin{lemma}}
\newcommand\EndLemma{\end{lemma}}
\newcommand\Proposition{\begin{proposition}}
\newcommand\EndProposition{\end{proposition}}
\newcommand\Corollary{\begin{corollary}}
\newcommand\EndCorollary{\end{corollary}}
\newcommand\Claim{\begin{claim}}
\newcommand\EndClaim{\end{claim}}
\newcommand\Proof{\begin{proof}}
\newcommand\EndProof{\end{proof}}
\newcommand\Equation{\begin{equation}}
\newcommand\EndEquation{\end{equation}}
\newcommand\NoProof{{\hfill$\square$}}
\newcommand\Bullets{\begin{itemize}}
\newcommand\EndBullets{\end{itemize}}

\newcommand\magentapacking{discrete set}
\newcommand\net{discrete set}

\newcommand\calT{{\mathcal T}}

\newcommand\fraka{\mathfrak a} 

\newcommand\frakc{\mathfrak c}

\renewcommand\epsilon{\varepsilon}

\newcommand\RR{{\bf R}}

\newcommand\redG{G}
\newcommand\totgeo{totally geodesic}

\newcommand\redS{S}
\newcommand\redGST{G^{S,\redcalT}}
\newcommand\reds{s}
\newcommand\redGndd{G _{\rm ndd}}
\newcommand\redGamma{\Pi}
\newcommand\wasbarX{A}
\newcommand\wasbeta{h}
\newcommand\wasbarY{B}

\newcommand\XPS{\calX_{P,\scrS}}

\newcommand\dotsystem{dot system}
\newcommand\adapted{adapted}
\newcommand\good{good}
\newcommand\veccale{{\vec\cale}}

\newcommand\inter{\mathop{\rm int}}

\newcommand\Mthick{M_{\rm thick}}
\newcommand\Mthin{M_{\rm thin}}
\newcommand\redasubs{\cala_S}

\newcommand\bareta{\overline{\eta}}
\newcommand\bartheta{\overline{\theta}}

\newcommand\Hstar{H^\ast}
\newcommand\Ustar{u^\ast}
\newcommand\init{\mathfrak{i}}
\newcommand\term{\mathfrak{t}}
\newcommand\Hterm{H_1}
\newcommand\Uterm{u_1}
\newcommand\Cterm{U_1}
\newcommand\Hinit{H_0}
\newcommand\Uinit{u_0}
\newcommand\Cinit{U_0}
\newcommand\wasU{u}
\newcommand\wouldbeV{v}
\newcommand\reddenedetanought{\eta_0}
\newcommand\reddenedetaone{\eta_1}
\newcommand\wasC{U}

\newcommand\ttau{\widetilde\tau}
\newcommand\tH{\widetilde H}

\newcommand\tU{{\widetilde U}}

\newcommand\maybecalD{\calD}
\newcommand\willbep{p}
\newcommand\tryepsilon{\epsilon}

\newcommand\willbeS{S}

\newcommand\willbefinitevolume{finite-volume}
\newcommand\willbesemifree{semifree}

\newcommand\willbehatp{{\widehat p}}
\newcommand\hatp{{\widehat p}}
\newcommand\willbepolytope{polyhedron}

\newcommand\willbepolytopal{polyhedral}

\newcommand\willbepolyhedra{polyhedra}

\newcommand\maybeepsilon{\epsilon}
\newcommand\maybeb{b}

\newcommand\tS{{\widetilde S}}

\newcommand\talpha{{\widetilde\alpha}}
\newcommand\txi{{\widetilde\xi}}
\newcommand\tbeta{{\widetilde\beta}}

\newcommand\calc{{\mathcal C}}
\newcommand\calD{{\mathcal D}}

\newcommand\cala{{\mathcal A}}

\newcommand\calF{{\mathcal F}}
\newcommand\calf{{\mathcal F}}

\newcommand\cale{{\mathcal E}}

\newcommand\scrGndd{{\scrG}_{\rm ndd}}

\newcommand\redGstar{\redGndd}
\newcommand\redscrPST{\scrP_{S,\redcalT}}

\newcommand\calenl{{\mathcal E}_{\rm nl}}
\newcommand\redcalestar{\cale_{\rm ndd}}

\newcommand\realcalX{{\mathcal X}}
\newcommand\realcalH{{\mathcal H}}
\newcommand\calX{X}

\newcommand\baralpha{\overline{\alpha}}

\newcommand\atwoS{{\cala^2_S}}
\newcommand\athreeS{{\cala^3_S}}

\newcommand\QQ{{\mathbb Q }}
\newcommand\HH{{\bf H}}
\newcommand\EE{{\bf E}}

\newcommand\scrD{{\mathscr D}}
\newcommand\scrG{{\mathcal G}}

\newcommand\scrP{{\mathscr P}}

\newcommand\scrS{{\mathscr S}}

\newcommand\scrV{{\mathscr V}}
\newcommand\scrL{{\mathscr L}}

\newcommand\dist{\mathop{\rm dist}}

\newcommand\nbhd{\mathop{\rm nbhd}}

\newcommand\density{{\rm density}}

\newcommand\vol{\mathop{\rm vol}}

\newcommand\length{\mathop{{\rm length}}}

\newcommand\gpc{generalized polyhedral complex}

\newcommand\rank{\mathop{{\rm rank}}}

\newcommand\image{\mathop{{\rm Im}}}

\newcommand\isomplus{\mathop{{\rm Isom}_+}}
\newcommand\isom{\mathop{{\rm Isom}}}

\newcommand\FF{{\bf F}}

\newcommand\calM{{\mathcal M}}
\newcommand\calH{{\mathscr H}}
\newcommand\redcalM{\mathcal M}
\newcommand\redcalT{\mathcal T}
 
\newcommand\scrGST{\scrG^{S,\redcalT }}
\newcommand\shortone{{\mathfrak s}_1}

\newcommand\mayber{r}

\newcommand\redsubone{_1}
\newcommand\redsubnought{_0}
\DeclareFontFamily{U}{rcjhbltx}{}
\DeclareFontShape{U}{rcjhbltx}{m}{n}{<->rcjhbltx}{}
\DeclareSymbolFont{hebrewletters}{U}{rcjhbltx}{m}{n}

\let\aleph\relax\let\beth\relax
\let\gimel\relax\let\daleth\relax

\DeclareMathSymbol{\aleph}{\mathord}{hebrewletters}{39}
\DeclareMathSymbol{\beth}{\mathord}{hebrewletters}{98}
\DeclareMathSymbol{\gimel}{\mathord}{hebrewletters}{103}
\DeclareMathSymbol{\daleth}{\mathord}{hebrewletters}{100}

\DeclareMathSymbol{\lamed}{\mathord}{hebrewletters}{108}
\DeclareMathSymbol{\mem}{\mathord}{hebrewletters}{109}
\DeclareMathSymbol{\ayin}{\mathord}{hebrewletters}{96}
\DeclareMathSymbol{\tsadi}{\mathord}{hebrewletters}{118}
\DeclareMathSymbol{\qof}{\mathord}{hebrewletters}{114}
\DeclareMathSymbol{\shin}{\mathord}{hebrewletters}{152}

\usepackage[utf8]{inputenc}

\begin{document}

\title{The ratio of homology rank to hyperbolic volume, II}

\author{Rosemary K. Guzman}
\address{Department of Mathematics
\\
University of Illinois\\
1409 W. Green St.\\
Urbana, IL 61801}
\email{rguzma1@illinois.edu}

\author{Peter B. Shalen}
\address{Department of Mathematics, Statistics, and Computer Science
(M/C 249)\\
University of Illinois at Chicago\\
851 S. Morgan St.\\
Chicago, IL 60607-7045}
\email{petershalen@gmail.com}

\maketitle

\begin{abstract}
Under  mild topological restrictions, we obtain new linear upper
bounds for the dimension of the mod $p$ homology (for any prime $p$)
of a finite-volume orientable hyperbolic $3$-manifold $M$ in terms of
its volume. A surprising feature of the arguments in the paper is that they
require an application of the Four Color Theorem.

If $M$ is closed, and either 
(a) $\pi_1(M)$ has no subgroup isomorphic to the fundamental group
  of a closed, orientable  surface of genus $2$, $3$ or $4$, or (b) $p=2$, and $M$ contains no (embedded, two-sided) incompressible
  surface of genus $2$, $3$ or $4$, then
$\dim H_1(M;\FF_p)< 157.763\cdot\vol (M)$. If $M$ has one or more
cusps, we get a very similar bound assuming that $\pi_1(M)$ has no subgroup isomorphic to the fundamental group
  of a closed, orientable  surface of genus $g$ for
  $g=2,\ldots,8$. These results should be compared with those of our
  previous paper ``The ratio of homology rank to hyperbolic volume,
  I,'' in which we obtained a bound with a coefficient in the range of
  $168$ instead of $158$, without a restriction on surface subgroups
  or incompressible surfaces. In a future paper, using a much more
  involved argument, we expect to obtain bounds close to those given by the
  present paper without such a restriction.

The arguments also give new linear upper
bounds (with constant terms) for the rank of $\pi_1(M)$ in terms of
$\vol M$, assuming that either $\pi_1(M)$ is $9$-free, or $M$ is
closed and $\pi_1(M)$ is $5$-free.
\end{abstract}

\section{Introduction}

It is a standard consequence of the Margulis Lemma
that there is a constant $\lambda>0$ such that for every
\willbefinitevolume
\ orientable
hyperbolic $3$-manifold and every prime $p$ we have
\Equation\label{in general}
\dim H_1(M;\FF_p)\le\lambda\cdot\vol M
\EndEquation
(where $\vol$ denotes hyperbolic volume).
In \cite[Proposition 2.2]{alm}, Agol, Leininger and Margalit 
established this inequality with a relatively good value of $\lambda$,
by using the
so-called $\log3$ Theorem---a special case of the ``$\log(2k-1)$
Theorem,'' which was put in final form in \cite{acs-surgery}, building on
preliminary versions in \cite{paradoxical} and \cite{accs}. (Theorem
\ref{another log(2k-1)} of this paper is a version of the $\log(2k-1)$
Theorem.)  The use of
the $\log(2k-1)$ Theorem in \cite{alm} allows one to improve on a more
naive value of $\lambda$ by orders of magnitude.

In \cite{ratioI} it was shown that the inequality (\ref{in general}) holds with a
significantly smaller value of $\lambda$ 
than the one given in
\cite{alm}. 
In this paper we introduce a
surprising new idea that gives  further improvements on the value of
$\lambda$ under  mild topological hypotheses on the manifold $M$.
In a future paper we will show that when $M$ is closed and the prime $p$ is equal to $2$, a similar improvement on the
result of \cite{ratioI} can be obtained without any additional
topological hypothesis, but the arguments needed to prove this are
much more involved than those of the present paper.

The surprising feature of the arguments in this paper is that they
require an application of the Four Color Theorem. This further
illustrates a pattern that was already seen in the papers
\cite{acs-surgery}, \cite{acs-singular}, \cite{singular-two},
\cite{fourfree}, and \cite{kfree-volume},
in which a broad range of techniques and results from pure mathematics
were
brought
to bear on the problem of bounding
$\dim H_1(M;\FF_p)$ in terms of a given bound on $\vol M$, where $M$
is a
\willbefinitevolume
\ orientable
hyperbolic $3$-manifold and $p$ is a prime. 
These papers invoke the $\log(2k-1)$ Theorem, which itself is proved
using a Banach-Tarski-style
decomposition of the Patterson-Sullivan measure for a free Kleinian
group 
the Marden tameness
conjecture, proved by Agol \cite{agol-tameness} and Calegari-Gabai
\cite{cg}. Other results and techniques invoked in the papers \cite{acs-surgery}, \cite{acs-singular}, \cite{singular-two},
\cite{fourfree}, and \cite{kfree-volume} include
Perelman's estimates for Hamilton's Ricci flow, via the
applications in \cite{ast}; the tower construction used by Shapiro and
Whitehead in their proof of Dehn's lemma (\cite{shapiro-whitehead});
the Borsuk nerve theorem; the topological theory of $3$-manifold group
actions on trees;
and Fisher's inequality on block designs.

For the case of a closed manifold, the main result of this paper
relating homology to volume is the following theorem:
\begin{introclosedhomologytheorem}
Let $M$ be a closed, orientable hyperbolic
$3$-manifold, and let $p$ be  a prime. Suppose that either 
\begin{enumerate}[(a)]
\item $\pi_1(M)$ has no subgroup isomorphic to the fundamental group
  of a closed, orientable  surface of genus $2$, $3$ or $4$; or
\item $p=2$, and $M$ contains no (embedded, two-sided) incompressible
  surface of genus $2$, $3$ or $4$. 
\end{enumerate}
Then 
$$\dim H_1(M;\FF_p)< 157.763\cdot\vol (M).$$
\end{introclosedhomologytheorem}

For the case where the manifold $M$ is not assumed to be closed, we
will prove the following result:

\begin{intromaybecuspstheorem}
Let $M$ be a finite-volume, orientable hyperbolic
$3$-manifold, and let $p$ be  a prime. Suppose that for
$g=2,3,\ldots,8$, the group $\pi_1(M)$ has no subgroup isomorphic to the fundamental group
  of a closed, orientable  surface of genus $g$.
Then 
$$\dim H_1(M;\FF_p)< 
158.12\cdot\vol( M).$$
\end{intromaybecuspstheorem}

Theorems \ref{closed homology} and \ref{maybe cusps} are deduced,
via results from \cite{accs},
\cite{singular-two}, and \cite{milley}.
from
theorems relating volume to the fundamental group. A group $\Pi$ is
said to be {\it $k$-free} for a given positive integer $k$ if each 
subgroup of $\Pi$ having rank at most $k$ is free. A group $\Pi$ is
said to be {\it $k$-semifree} for a given  $k$ if each 
subgroup of $\Pi$ having rank at most $k$ is a free product of free
abelian groups. We will prove the following two results:

\begin{intromainertheorem}
Let $M$ be a \willbefinitevolume\  orientable hyperbolic $3$-manifold
such that $\pi_1(M)$ is $5$-free. 
Then we have
$$\rank\pi_1(M)<1+
157.497\cdot
\vol(M).$$
\end{intromainertheorem}

\begin{introsemimainertheorem}
Let $M$ be a \willbefinitevolume\  orientable hyperbolic $3$-manifold
such that $\pi_1(M)$ is $9$-semifree. 
Then we have
$$\rank\pi_1(M)<1+
157.766
\cdot
\vol(M)
.$$
\end{introsemimainertheorem}

Theorems \ref{closed homology} and \ref{maybe cusps} are partial
improvements on Theorem 5.4 of \cite{ratioI}, while Theorems 
\ref{mainer} and \ref {semimainer} are partial improvements on  
 Proposition 5.2 of \cite{ratioI}. Theorem 5.4 and Proposition 5.2 of
 \cite{ratioI} have weaker hypotheses than their counterparts in this
 paper, in that they do not require restrictions on incompressible
 surfaces in the manifold $M$ or surface subgroups of $\pi_1(M)$. On the other hand, the results in \cite{ratioI} have weaker
 conclusions, in that the coefficients in the inequalities relating
 volume to the dimension of $H_1(M;\FF_p)$ or the rank of $\pi_1(M)$
 are in the range of $168$ rather than $158$.
The reason for the improvement by a factor of about
$15/16$ over the constants in \cite{ratioI} will emerge in the sketch
given below of the proofs of Theorems \ref{mainer} and
\ref{semimainer}, and will be seen to involve the Four Color Theorem.

As context for sketching the new ideas in the proofs of
Theorems \ref{mainer} and \ref{semimainer},
we briefly review the methods of \cite{alm} and \cite{ratioI}.

The hypothesis that
the \willbefinitevolume\ orientable $3$-manifold 
$M$ has a $2$-\willbesemifree\ fundamental group implies, via the
$\log3$ Theorem, that $\tryepsilon:=\log3$ is a Margulis
number for $M$. One fixes a
finite subset $S$ of $M$ which is a maximal
$\tryepsilon$-\magentapacking\ 
contained in the
$\tryepsilon$-thick part of $M$. 
By a
Voronoi-Dirichlet construction, one associates with the
set $S$ a ``cell complex'' 
 $K_S$ whose underlying space is $M$. (The notions of Margulis number, 
$\epsilon$-thick part and $\epsilon$-discrete set, and the
construction and properties of $K_S$, are developed in
some detail in \cite{ratioI} and are briefly reviewed in
Section 4 of the present paper.
  When $M$ is closed, $K_S$ is a
 certain kind of CW
 complex.) Each (open)
$3$-cell of $K_S$ is the homeomorphic image  under a locally isometric
covering map $q:\HH^3\to M$ of the interior of a
``Voronoi region,''
which is a convex \willbepolytope\  $X$ associated with a point $P$ of
$q^{-1}(S)$, and having $P$ as an interior point. We think of $P$ as
the ``center'' of the Voronoi region $X$. The Voronoi regions and
their faces form a {\it polyhedral complex,} in a sense that is defined in 
\cite{ratioI}
and reviewed in
Section 3 of the present paper.
This polyhedral complex has underlying space $\HH^3$, and will be denoted by $\realcalX$
in the present sketch; the cells of $K_S$ are the images
under $q$ of the interiors of the polyhedra belonging to $\realcalX$.

Now
$\pi_1(M)$ is carried by 
a connected graph $\scrG$ 
which has
$S$ as its vertex set, and has one edge ``dual'' to each $2$-cell of $K_S$
meeting the 
$\epsilon$-thick 
part of $M$.
To bound the rank of $\pi_1(M)$ it suffices
to bound the first betti number of $\scrG$; this in turn can be done
by bounding the number of vertices of $\scrG$ and the valence of an
arbitrary vertex of $\scrG$. A bound on the number of vertices is
obtained in \cite{alm} by using that $\epsilon$ is a Margulis number
for $M$, and is improved in \cite{ratioI} by a simple
application of B\"or\"oczky's density
bound for sphere-packings in hyperbolic space \cite{boroczky}. The
main new idea in \cite{ratioI} is a refinement, involving novel,
elementary geometric arguments, of the rather straightforward method
used in \cite{alm} to bound the valence of a vertex in $\scrG$.

The arguments of the present paper build on the bounds given in
\cite{ratioI} for the first betti number of $\scrG$. 
The basic
strategy for improving the bounds for the rank of $\pi_1(M)$ given in
\cite{ratioI} is to try to construct a connected subgraph of $\scrG$ which has
significantly smaller betti number than $\scrG$ but still carries the
fundamental group of $M$. While one does not always succeed in
constructing such a subgraph, one can come close enough for the
purpose of proving Theorems \ref{mainer} and \ref {semimainer}, as the
following sketch indicates.

The context for the arguments of this paper involves the use of a
Margulis number $\epsilon$ which is slightly less than $\log3$. Both
the number $\epsilon$ and the maximal $\epsilon$-discrete set $S$ in the $\epsilon$-thick
part of $M$ are taken to be generic in a suitable sense. One of the
genericity properties of $S$ is that the polyhedral complex
$\realcalX$ which it defines is {\it weakly simple} in the sense that
every $1$-dimensional polyhedron in $\realcalX$ is a face of exactly
three $3$-dimensional polyhedra in $\realcalX$; the proof that weak
simplicity is indeed a generic property depends on work by Kapovich
\cite{kapovich}.

Having chosen a generic $S$, we use the Four Color Theorem to color
the two-dimensional faces of each three-dimensional polyhedron
$X\in\realcalX$, using only four colors, in such a way that any two two-dimensional faces of
$X$ which share a one-dimensional face have distinct colors. This
coloring is done in a $\Gamma$-invariant fashion, where $\Gamma$ is
the group of deck transformations of the covering $q:\HH^3\to M$. 

Now consider any oriented edge $\eta$ of the dual graph $\scrG$ to
$K_S$. If $\alpha$ is a path in $M$ realizing $\eta$, a lift $\talpha$
of $\alpha$ to $\HH^3$ has its initial point in a $3$-dimensional
polyhedron $X\in\realcalX$, and crosses the boundary of $X$ in the
interior of a two-dimensional face of $X$, which has now been given a
color. In this way, each oriented edge $\eta$ of the dual graph $\scrG$ 
determines a color.  

By a {\it color assignment} for $S$ we mean a function that assigns a
``distinguished'' color to each vertex of $\scrG$. Given a color
assignment, we shall say that an oriented edge $\eta$ of $\scrG$ is
{\it distinguished} if the color determined by $\eta$ is the
distinguished color for its initial vertex. An edge of $\scrG$ will be
termed {\it doubly distinguished} if both its orientations are
distinguished. 

If we make the simplifying assumption that the graph $\scrG$ has no
loops, then a counting argument shows that there is a color assignment
for $S$ with the property that at least one-sixteenth of the edges of
$\scrG$ are doubly distinguished. But one can use the genericity
properties of $S$, including weak simplicity, to show that the
subgraph $\scrGndd$ of $\scrG$, consisting of all vertices and all
non-doubly-distinguished edges of $\scrG$, is connected and carries the fundamental
group of $M$. This leads to an upper bound for the rank of $\pi_1(M)$
which is approximately fifteen-sixteenths of the upper bound for the
first betti number of $\scrG$ established in \cite{ratioI}.

In the general case, when $\scrG$ may contain loops, we consider the
subgraph $\scrG^\dagger$ of $\scrG$ consisting of all vertices and all
non-doubly-distinguished edges of $\scrG$ which are not loops. We also
consider, for each vertex $v$ of $\scrG$, the subgroup $B_v$ of $\pi_1(M,v)$
generated by elements that are defined by loops of the graph $\scrG$
based at $v$. It can
be shown that $\scrG^\dagger$ is connected and contains at most fifteen-sixteenths of
the edges of $\scrG$, and that $\pi_1(M)$ is generated by the image of
the inclusion $\pi_1(\scrG^\dagger)\to\pi_1(M)$, together with
isomorphic copies of the groups $B_v$ as $v$ ranges over the vertices
of $\scrG$. The additional information needed to deduce the conclusions
of Theorems \ref{mainer} and \ref {semimainer} is a bound on the rank
of $B_v$ for each vertex $v$ of $\scrG$. In the context of the proof
of Theorem
\ref{mainer}, we show that each loop of $\scrG$
defines a closed path in $M$ of length at most $2\epsilon<\log9$, and
using the case $k=5$ of the $\log(2k-1)$ Theorem, 
together with the
hypothesis that $\pi_1(M)$ is $5$-free, to conclude that each $B_v$
has rank at most $4$. A similar argument, in the context of the proof
of Theorem \ref{semimainer}, gives a bound of $8$ for the rank of each $B_v$.

The details of the arguments sketched above are
assembled in Sections 6 and 7, after developing the ingredients in
earlier sections. Section 2 supplements the material about convex
polyhedra that was presented in \cite{ratioI} by establishing some
additional facts that are needed in the present paper. In a similar
way, Section 4 supplements the material about Margulis numbers,
$\epsilon$-discrete sets, Voronoi complexes and dual graphs that was
presented in  \cite{ratioI}. Section 3 contains the arguments that are
needed for showing that weak simplicity is a generic property for a
Voronoi complex. The facts about colorings and color assignments that
we need are proved in Section 5. (The main results of Section 5,
Propositions \ref{because four colors}, \ref{one-sixteenth} and
\ref{homotoping prop}, are quoted in the proof of a central result,
Proposition \ref{new what's new}.)

In a future paper we will establish a result which has almost as
strong a conclusion as the case $p=2$ of Theorem \ref{closed
  homology}, but does not require the assumption that there are no
closed incompressible surfaces of genus $2$, $3$ or $4$. The proof of
this result, which is much more involved than that of Theorem \ref{closed
  homology}, combines the methods of the present paper with the
results of \cite{ast} and \cite{miyamoto}, and desingularization
techniques similar to those of \cite{acs-singular} and \cite{singular-two}.

We summarize here some conventions that will be used in the body of
the paper.

A {\it path} in a space $X$ is a map $\alpha:[0,1]\to X$. A {\it
  closed path} in $X$ is a path $\alpha$ such that
$\alpha(0)=\alpha(1)$; the term {\it loop} is reserved for its
graph-theoretical meaning (see \ref{new graph stuff}). The composition
of two paths $\alpha$ and $\beta$, with $\alpha(1)=\beta(0)$, will be
denoted by $\alpha\star\beta$ (so that $\alpha\star\beta\,(t)$ is equal
to $\alpha(2t)$ if $0\le t\le 1/2$ and to $\beta(2t-1)$ if $1/2\le
t\le 1$). If $\alpha$ is a path, $\baralpha$ will denote the path defined by
$\baralpha(t)=\alpha(1-t)$. 

The symbol  ``$\dist$'' will denote the distance in a metric
space when it is clear which metric space
is involved.
We denote the group of isometries of a
metric space $X$ by $\isom(X)$; if $X$ is a connected, orientable
Riemannian manifold, then $\isomplus(X)$ will denote the group of
orientation-preserving isometries of $X$. 

If $X$ is a group, we write $Y\le X$ to mean that $Y$ is a subgroup of $X$.

\section{Convex polyhedra}

\Number\label{in the beginning}
In \cite[Subsection 2.1]{ratioI}, a {\it convex set} in
$\EE^n$ or
 $\HH^n$ is defined to be a non-empty set which is the intersection of a family
 of half-spaces, and a {\it convex polyhedron} is defined to be a non-empty set which is the intersection of a family
 of half-spaces whose bounding hyperplanes form a locally
 finite family.
The {\it interior} $\inter X$ of a convex set $X$ is defined to be its
topological interior relative to the smallest \totgeo\ subspace $W$ of 
$\EE^n$ or
$\HH^n$
containing $X$, and the {\it dimension} $\dim X$ of $X$ is defined to be the
dimension of $W$.
A {\it support hyperplane} for a closed convex set  $X$ is
a hyperplane  which has non-empty intersection with $X$,
and is the boundary of a closed half-space containing $X$. 
A {\it
  face} of $X$ is defined to be
a subset $F$ of $\EE^n$ or $\HH^n$ 
which 
either is
  equal to $X$, 
or
is the intersection of $X$ with a support hyperplane for $X$ (and in  the
latter case  $F$ is called a {\it proper} face).

As in \cite[Subsection 2.1]{ratioI}, we shall say that 
a family $\calH$ of closed
half-spaces in $\HH^n$ is {\it irredundant} if $\calH$ has no proper subfamily $\calH'$ such that
$\bigcap_{\realcalH\in\calH'}\realcalH
=\bigcap_{\realcalH\in\calH}\realcalH$. 

As in \cite[Subsection 2.1]{ratioI}, we define a
{\it facet} of a convex polyhedron $X$ to be a maximal proper face of
$X$. According to \cite[Assertion 2.1.9]{ratioI}, the facets of $X$
are precisely its codimension-$1$ faces.
\EndNumber

\Number\label{cardinal}
In \cite[Subsection 2.1]{ratioI}
we listed a number of 
facts about 
closed
convex sets and convex  \willbepolyhedra\ in $\HH^n$ that are the
counterparts of 
well-known facts about 
closed
convex  sets and convex \willbepolyhedra\  in
$\EE^n$.  
Here we will point out a few 
additional such facts which will be needed in the present paper, and give brief hints
about how to prove them. 

\Claim\label{finite hull}
Let $T$ be a finite subset of $\HH^n$. Then the convex
hull of $T$, i.e. the smallest convex subset of $\HH^n$ containing
$T$, is a compact convex hyperbolic polyhedron.
\EndClaim

The Euclidean analogue of \ref{finite hull} is included in
\cite[Section 3.1, Assertion 1]{Grunbaum}. Indeed, according to the
latter assertion, every ``polytope'' in $\EE^n$, in the sense defined in
\cite[Section 3.1]{Grunbaum}, is a finite intersection of closed
half-spaces, and hence in particular it is a convex polyhedron in the sense
defined above; and it is pointed out on page 31 of \cite{Grunbaum}
that ``polytopes'' may be characterized as convex hulls of finite
sets. Using the Beltrami model, it is easy to deduce \ref{finite hull}
from its  Euclidean analogue.

\Claim\label{awk}
Let $X$ be a convex hyperbolic polyhedron in $\HH^n$ which is given in
the form
$X=\bigcap_{\realcalH\in\calH}\realcalH$, where $\calH$ is a family of closed
half-spaces in $\HH^n$ whose bounding hyperplanes form a locally
finite family. Suppose that $\calH$ is irredundant. Then the facets of $X$ are precisely the sets of the form $X\cap\partial \realcalH$ for
$\realcalH\in\calH$.
\EndClaim

In the special case where the family $\calH$ is finite, the Euclidean
analogue of \ref{awk} is included in \cite[Section 2.6, Assertions 2
and 3]{Grunbaum}. The argument for this special case goes through
without change in the hyperbolic context, and it is a straightforward exercise to extend the proof from  this special case
to the general  case.

If a convex hyperbolic polyhedron $X$  in $\HH^n$  is given in
the form
$X=\bigcap_{\realcalH\in\calH}\realcalH$, where $\calH$ is a family of closed
half-spaces in $\HH^n$ whose bounding hyperplanes form a locally
finite family, then according  \cite[Assertion 2.1.7]{ratioI}, $\calH$
has an irredundant subfamily $\calH_0$ such that
$X=\bigcap_{\realcalH\in\calH_0}\realcalH$. Combining this with
\ref{awk}, we immediately deduce:

\Claim\label{kwa}
Let $X$ be a convex hyperbolic polyhedron in $\HH^n$ which is given in
the form
$X=\bigcap_{\realcalH\in\calH}\realcalH$, where $\calH$ is a family of closed
half-spaces in $\HH^n$ whose bounding hyperplanes form a locally
finite family. Then every facet of $X$ has the form
$X\cap\partial \realcalH$ for some
$\realcalH\in\calH$.
\EndClaim

The following fact is an easy consequence of \cite[Assertion
2.1.6]{ratioI}, which asserts that the faces of a convex \willbepolytope\ form
a locally finite family of sets.

\Claim\label{in a facet}
Every proper face of a convex polyhedron $X$ is contained in a facet of $X$.
\EndClaim

\EndNumber

\Proposition\label{Toss in vertices}
Let $C$ be a compact subset of $\HH^n$ for a given $n\ge2$. Then there
   is a compact convex
hyperbolic polyhedron $G\subset\HH^n$ with $\dim C=n$ and $C\subset \inter G$.
\EndProposition

\Proof
Every point of $\HH^n$ lies in the interior of a hyperbolic
$n$-simplex. Since $C$ is compact, there is a
finite collection $\scrD$ of hyperbolic $n$-simplices such that
$C\subset\bigcup_{\Delta\in\scrD}\inter\Delta$. We may take $\scrD$ to
be non-empty. The vertices of the
simplices in $\scrD$ form a non-empty finite subset $T$ of $\HH^n$. 
According to \ref{finite hull}, the convex
hull of $T$ is a compact convex hyperbolic polyhedron $G$; and since
$G$ contains at least one $n$-simplex, we have $\dim G=n$. We have
$C\subset\inter G$.
\EndProof

\Proposition\label{facets in facets}
Let $X_1$ and $X_2$ be convex polyhedra in $\HH^n$ for a given
$n\ge2$. Suppose that $X_1\cap X_2\ne\emptyset$. Then $X_1\cap X_2$ is a
convex polyhedron, and each proper face of $X_1\cap X_2$ is contained  in a
proper face of one of the $X_i$.
\EndProposition

\Proof
For $i=1,2$, let us write $X_i=\bigcap_{\realcalH\in\calH_i}\realcalH$, where
$\calH_i$ is a family of half-spaces whose bounding hyperplanes form a locally
finite family. It follows from \cite[Assertion 2.1.7]{ratioI} that we
may take each of the families $\calH_1,\calH_2$ to be
irredundant.
Now $\calH:=\calH_1\cup\calH_2$ is
again a family of half-spaces whose bounding hyperplanes form a locally
finite family, and $X:=X_1\cap X_2=\bigcap_{\realcalH\in\calH}\realcalH$; this shows
that $X$ is a convex polyhedron, which is the first assertion of the
proposition. To prove the second assertion, in view of \ref{in a
  facet}, it suffices to show
that every facet of $X$ is
contained in a facet of one of the $X_i$.  If $F$ is a facet of $X$,
then by \ref{kwa} we may write $F=X\cap\partial \realcalH_0$ for some hyperplane
$\realcalH_0\in\calH$. Then for some $i_0\in\{1,2\}$ we have $\realcalH_0\in\calH_{i_0}$. We
have $G:= X_{i_0}\cap\partial \realcalH_0 \supset F$; and in view of the irredundancy
of $\calH_{i_0}$, it follows from \ref{awk} that  $G$ is a facet of
$X_{i_0}$.
\EndProof


\Lemma\label{cincinnati}
If $C$ is a compact subset of a
convex hyperbolic polyhedron $F$, then there is a compact convex
hyperbolic polyhedron $F_0\subset F$, having the same
dimension as $F$, such that $C\subset F_0$. 
Furthermore, each
proper face of $F_0$ is either a proper face of $F$ or is disjoint from $C$. 
\EndLemma

\Proof
After possibly replacing $C$ by a larger compact subset of $F$, we may
assume $C\cap\inter F\ne\emptyset$. According to Proposition \ref{Toss in vertices}, there
   is a compact convex
hyperbolic polyhedron $G\subset\HH^n$ with $\dim G=n$ and $C\subset
\inter G$. According to Proposition \ref{facets in facets}, $F_0=F\cap
G$ is either a convex polyhedron or the empty set; and since $G$ is compact, $F_0$ is also
compact. We have $\inter F\cap\inter G\supset\inter F\cap
C\ne\emptyset$. Hence $F_0$ is a convex polyhedron, and contains a
nonempty open subset of the smallest \totgeo\ subspace of $\HH^n$
containing $F$. This implies that $\dim F_0=\dim F$. 
Since $C\subset F$ and
$C\subset\inter G$, we have $C\subset F_0$.
If $E$ is a proper face of $F_0$, then by Proposition \ref{facets in facets},
$E$ is a proper face of either $F$ or $G$; and if $E$ is a proper face of $G$, then
since $C\subset\inter G$, we have $E\cap C=\emptyset$.
\EndProof

\section{Weakly simple Voronoi complexes}

\Number\label{complex review}
We recall from \cite{ratioI} that a {\it \willbepolytopal\ complex} in
$\HH^n$ (for any $n\ge2$) is defined to be a locally finite
collection $\realcalX$ of convex \willbepolyhedra\ in $\HH^n$ such
that (i) for each $X\in\realcalX$, all the faces of $X$ belong to
$\realcalX$, and (ii) for all $X,Y\in\realcalX$, the set $X\cap Y$
either is a common face of $X$ and $Y$ or is empty. 
As in \cite{ratioI} we denote by
$|\realcalX|$ the union of all the convex \willbepolyhedra\ that
belong to a given polyhedral complex $\realcalX$.
\EndNumber

\Lemma\label{before about 3}
Suppose that $\realcalX$
is a polyhedral complex such that $|\realcalX|=\HH^3$. Then for every three-dimensional polyhedron
$X\in\realcalX$ and every one-dimensional face $L$ of  $X$, there are
three-dimensional polyhedra $X_1,X_2\in\realcalX$, distinct from $X$ and
from each other, such that 
\begin{itemize}
\item
for $i=1,2$, 
the set $X\cap X_i$ is a common
two-dimensional face of $X$ and $X_i$, 
\item
for $i=1,2$, the polyhedron $L$ is a face of $X\cap X_i$,
and
\item
$X\cap X_1\cap X_2=L$.
\end{itemize}
\EndLemma

\Proof
It follows from \cite[Assertion 2.1.1]{ratioI} that $\partial X$ is a
$2$-manifold. Hence
there are two distinct two-dimensional faces $F_1$ and
$F_2$ of $X$ having $L$ as a face, and we have $F_1\cap F_2=L$.

Now let an index $i\in\{1,2\}$ be given. Choose an interior point $P$ of $F_i$. Since $\realcalX$ is
a polyhedral complex and $|\realcalX|=\HH^3$, some neighborhood of $P$
in $\HH^3$ is contained in a  union of polyhedra in $\realcalX$
having $F_i$ as a face. Since $F_i\cup X$ is not itself a neighborhood
of $P$, there must be a three-dimensional polyhedron in $\realcalX$
which is distinct from $X$ and has $F_i$ as a face. We choose such a
polyhedron and denote it by $X_i$. Now since $X$ and $X_i$ are
distinct three-dimensional polyhedra in the polyhedral complex $\realcalX$,
their intersection is a priori empty or a proper face of both $X$ and
$X_i$; since this intersection contains the common two-dimensional
face $F_i$ of $X$ and $X_i$, we must have $X\cap X_i=F_i$. Thus $L$ is
a face of $X\cap X_i$.

We have  $X\cap X_1\cap X_2=F_1\cap F_2=L$.

Since we have $X\cap X_i=F_i$ for $i=1,2$, and since $F_1\ne\ F_2$, we
must have $X_1\ne X_2$. All the asserted properties of $X_1$ and $X_2$
are now established.
\EndProof

\Definition\label{valence def}
If $\realcalX$ is a polyhedral complex such that $|\realcalX|$ is a
$3$-manifold without boundary, and $L\in\realcalX$ is a
one-dimensional polyhedron, we define the {\it valence} of $L$ to be
the number of three-dimensional polyhedra which are elements of
$\realcalX$ and have $L$ as a face.
\EndDefinition

\Proposition\label{about 3}
Suppose that $\realcalX$
is a polyhedral complex such that $|\realcalX|=\HH^3$. Then for every
one-dimensional polyhedron $L\in \realcalX$, the
  following conclusions hold.
\begin{enumerate}
\item  The valence of $L$ is at least $3$.
\item The intersection of all three-dimensional polyhedra having $L$ as a
face is $L$.
\item If $L$ has valence exactly $3$,
then for any two distinct three-dimensional polyhedra $X,X'$ having $L$ as a
face, the set $X\cap X'$ is a two-dimensional common face of $X$ and
$X'$, and has $L$ as a face.
\end{enumerate}
\EndProposition

\Proof
To prove the first two assertions, we begin by choosing an interior point $P$ of 
the one-dimensional polyhedron $L\in \realcalX$. Since $\realcalX$ is
a polyhedral complex and $|\realcalX|=\HH^3$, some neighborhood of $P$
in $\HH^3$ is contained in a finite union of polyhedra in $\realcalX$
having $L$ as a face. 
Since an open set in $\HH^3$ cannot be a finite
union of convex polyhedra of dimension at most $2$, the polyhedron $L$ is a face of some
three-dimensional polyhedron $X_0\in\realcalX$. Lemma
\ref {before about 3}, applied with $X=X_0$, then gives, in particular, three-dimensional polyhedra $X_1,X_2\in\realcalX$, distinct from $X_0$ and
from each other, having $L$ as a face; 
this proves Conclusion (1).  Lemma
\ref {before about 3} also provides the information that $X_0\cap
X_1\cap X_2=L$, which establishes Conclusion (2).

To prove (3), suppose that the one-dimensional polyhedron $L\in \realcalX$ has
valence exactly $3$, and that  $X$ and $X'$
are distinct three-dimensional polyhedra having $L$ as a
face. Lemma \ref {before about 3} gives three-dimensional polyhedra $X_1,X_2\in\realcalX$, distinct from $X$ and
from each other, such that for $i=1,2$ the set $X\cap X_i$ is a common
two-dimensional face of $X$ and $X_i$, and $L$ is a face of $X\cap
X_i$. Since $L$ has valence $3$, the three distinct polyhedra $X$,
$X_1$ and $X_2$ are the only three-dimensional polyhedra in $\realcalX$
having $L$ as a face. Since $X'\ne X$ we must have $X'=X_{i_0}$ for some
$i_0\in\{1,2\}$. In particular,
$X\cap X'=X\cap X_{i_0}$ is a common
two-dimensional face of $X$ and $X_{i_0}$, i.e. of $X$ and $X'$; and $L$ is a face of $X\cap
X_{i_0}=X\cap X'$.
\EndProof

\DefinitionRemark\label{weakly simple def.}
Let $\realcalX$ be a polyhedral complex with $|\realcalX|=\HH^3$. We
will say that $\realcalX$ is {\it weakly simple} if every one-dimensional
polyhedron in $\realcalX$ has valence $3$.

It follows from Conclusion (1) of Proposition
\ref{about 3} that  a polyhedral complex $\realcalX$ with
$|\realcalX|=\HH^3$  is  weakly simple if and only if  every one-dimensional
polyhedron in $\realcalX$ has valence {\it at most} $3$.
\EndDefinitionRemark

\Number\label{voronoi review}
Let $\scrS$ be a locally finite (i.e. discrete and closed) subset of
$\HH^n$. As in \cite[Subsection 2.3]{ratioI}, for each point $P\in \scrS$, we set
$$\XPS =\{W\in\HH^n:\dist(W,P)\le\dist(W,Q)\text{ for every
}Q\in \scrS\}.$$
As in \cite{ratioI}, we shall often write $\calX_P$ in place of $\XPS $ in situations where
it is understood which  set $\scrS$ is involved.
We recall from \cite{ratioI} that a set having the form $\XPS $ for some $P\in \scrS$ will be called a
{\it Voronoi region} for the locally finite set $\scrS$.

According to \cite[Proposition 2.4]{ratioI}, each Voronoi region is a
convex polyhedron. As in \cite[Subsection 2.5]{ratioI}, we denote by
$\realcalX_\scrS$ the set of all  faces of
Voronoi regions for $\scrS$. According to \cite[Proposition
2.6]{ratioI}, $\realcalX_\scrS$ is a polyhedral complex.

One important special case of this construction is the one in which
$\scrS=\Gamma\cdot P$, where $P$ is a point of $\HH^3$ and
$\Gamma\le\isomplus(\HH^3)$ is a discrete group. The polyhedral
complex $\realcalX_\scrS$ is then called a Dirichlet complex. We will
use the following result, Theorem \ref{from kapovich}, concerning
Dirichlet complexes. This theorem is due to   M. Kapovich, and is
included in  Theorem 1.6 of his paper \cite{kapovich}. Indeed, 
\cite[Theorem 1.6]{kapovich} has exactly the same conclusion as
Theorem \ref{from kapovich}, and has a strictly weaker hypothesis:
instead of assuming that the discrete group $\Gamma$ is torsion-free,
Kapovich assumes only that it is ``of class $\mathcal K$,'' meaning
that its $2$-torsion satisfies a certain restriction. (Although the conclusion of \cite[Theorem 1.6]{kapovich}
is stated somewhat informally, the proof given in \cite{kapovich} uses
the language of weak simplicity, and gives exactly the conclusion stated
in Theorem \ref{from kapovich}
below.)
\EndNumber

\Theorem\label{from kapovich}
Let $\Gamma$ be a discrete, torsion-free subgroup of
$\isomplus(\HH^3)$. Then there is a dense $G_\delta$ set $\calT\subset\HH^n$ such that for
every point $P\in \calT$ the Dirichlet complex $\realcalX_{\Gamma\cdot P}$
is weakly simple. \NoProof
\EndTheorem

Theorem \ref{from kapovich} will play a role in the proof of the
following result:

\Proposition\label{when would it be}
Let $\Gamma$ be a discrete, torsion-free subgroup of
$\isomplus(\HH^3)$, and let $N$ be a positive integer. Then the space
$(\HH^3)^N$ contains a  dense $G_\delta$ set $Y$ such that, for each
element $(P_1,\ldots,P_N)$  of $Y$, the following conditions hold:
\begin{enumerate}
\item the orbits $\Gamma\cdot P_1,\ldots,\Gamma\cdot P_N$ are distinct; and
\item if we denote by $\scrS$ the locally finite family $\Gamma\cdot
  P_1\cup\cdots\cup\Gamma\cdot P_N$, the polyhedral complex
  $\realcalX_\scrS$ is weakly simple.
\end{enumerate}
\EndProposition

\Proof
Set $M=\HH^3/\Gamma$, and let $q:\HH^3\to M$ denote the quotient
map. Then $q$ defines a covering map $q^N:(\HH^3)^N\to M^N$. If we let $U\subset M^N$ denote the set of all $N$-tuples
$(p_1,\ldots,p_N)$ such that $p_i\ne\ p_j$ whenever $1\le i<j\le N$,
then $U$ is an open dense subset of $M^N$, and hence $\tU:=(q^N)^{-1}(U)$
is an open dense subset of $(\HH^3)^N$. Note that Condition (1) of the
statement holds
whenever $(P_1,\ldots,P_N)\in\tU$.

For each $(P_1,\ldots,P_N)\in(\HH^3)^N$ we set $\scrS(P_1,\ldots,P_N)=\Gamma\cdot
  P_1\cup\cdots\cup\Gamma\cdot P_N$.

Let $\Phi$ denote the set of all four-element subsets of the product
set $\{1,\ldots,N\}\times\Gamma$. For each
$(P_1,\ldots,P_N)\in(\HH^3)^N$, and each $\phi\in\Phi$, we set
$$E_\phi^{(P_1,\ldots,P_N)}=\bigcap_{(i,\gamma)\in\phi}X_{\gamma\cdot P_i}^{\scrS(P_1,\ldots,P_N)}.$$
For each $(i,\gamma)\in\phi$ the set $X_{\gamma\cdot P_i}^{\scrS(P_1,\ldots,P_N)}$ is a
Voronoi region for $\scrS(P_1,\ldots,P_N)$; hence
$E_\phi^{(P_1,\ldots,P_N)}$ is either the empty set or an element of
the polyhedral complex   $\realcalX_{\scrS(P_1,\ldots,P_N)}$.

For each $\phi\in\Phi$, we denote by $V_\phi$ the set of all points
$(P_1,\ldots,P_N)\in(\HH^3)^N$ such that $E_\phi^{(P_1,\ldots,P_N)}$
either is the empty set or consists of a single point. We claim:
\Claim\label{that's why}
For every $(P_1,\ldots,P_N)\in\bigcap_{\phi\in\Phi}V_\phi$, the
polyhedral complex $\realcalX_{\scrS(P_1,\ldots,P_N)}$ is weakly simple.
\EndClaim

To prove \ref{that's why}, let
$(P_1,\ldots,P_N)\in\bigcap_{\phi\in\Phi}V_\phi$ be given. If
$\realcalX_{\scrS(P_1,\ldots,P_N)}$ is not weakly simple, 
then by an observation made in \ref{weakly simple
  def.},
some one-dimensional polyhedron
$L\in\realcalX_{\scrS(P_1,\ldots,P_N)}$ has valence strictly greater
than $3$. We choose four distinct Voronoi
regions (i.e. three-dimensional polyhedra) in
$\realcalX_{\scrS(P_1,\ldots,P_N)}$, which
may 
be written in the form $X_{\gamma_j\cdot
  P_{i_j}}=X_{\gamma_j\cdot P_{i_j}}^{\scrS(P_1,\ldots,P_N)}$ for
$j=1,2,3,4$. Since the Voronoi regions $X_{\gamma_1\cdot
  P_{i_1}},\ldots, X_{\gamma_4\cdot P_{i_4}}$ are distinct, in
particular the pairs $(i_1,\gamma_1),\ldots, (i_4,\gamma_4)$ are
distinct; hence $\phi:=\{(i_1,\gamma_1),\ldots,
(i_4,\gamma_4)\}\in\Phi$, so that 
$(P_1,\ldots,P_N)\in V_\phi$.
By definition this means that the intersection 
$\bigcap_{i=1}^4 X_{\gamma_j\cdot
  P_{i_j}}$ is at most a single point, a contradiction since this
intersection contains $L$. Thus \ref{that's why} is proved.

We set 
$Y=\tU\cap\bigcap_{\phi\in\Phi}V_\phi$. Since $Y\subset\tU$, Condition
(1) of the conclusion of the proposition holds whenever
$(P_1,\ldots,P_N)\in Y$. It follows from \ref{that's why} that
Condition
(2) also holds whenever
$(P_1,\ldots,P_N)\in Y$. The rest of the proof will be devoted to
showing that $Y$ contains a dense $G_\delta$ set.

If $\phi$ is an element of $\Phi$, and $i$ is an index in
$\{1,\ldots,N\}$, we define the {\it multiplicity} of $i$ in $\phi$,
denoted $\mu(i,\phi)$, to be the number of ordered pairs in
$\phi\subset\{1,\ldots,N\}\times\Gamma$ whose first entry is $i$. We
then have 
$$\sum_{i=1}^N\mu(i,\phi)=4 \text{ for every }\phi\in\Phi.$$
When $4$ is expressed as a sum of non-negative integers, at least one
term is equal to $1$, $2$ or $4$. Hence if for each $m\ge0$ we set 
$$\Phi_i=\{\phi\in\Phi: \mu(i,\phi)=m\text{ for some }
  i\in\{1,\ldots,N\}\},$$
then 
\Equation\label{dragnet}
\Phi=\Phi_1\cup\Phi_2\cup\Phi_4.
\EndEquation
(The union on the right hand side of (\ref{dragnet}) is not a disjoint
union in general, because $\Phi_1$ and $\Phi_2$ are not disjoint.)

We claim:
\Claim\label{bridget}
The set $\bigcap_{\phi\in\Phi_4}V_\phi$ contains a dense $G_\delta$
set in
$(\HH^3)^N$.
\EndClaim

To prove \ref{bridget}, we first note that, since $\Gamma$ is discrete
and torsion-free, it follows from 
Theorem \ref{from kapovich}
that there is a dense $G_\delta$ set $\calT\subset\HH^3$ such that for
every point $P\in \calT$ the polyhedral complex $\realcalX_{\Gamma\cdot P}$
is weakly simple. Since  $\calT\subset\HH^3$ is a dense $G_\delta$ set, the
set $\calT^N\subset(\HH^3)^N$ is also a dense $G_\delta$ set. We will complete
the proof of \ref{bridget} by showing that $\calT^N\subset\bigcap_{\phi\in\Phi_4}V_\phi$.

Let $(P_1,\ldots,P_N)\in \calT^N$ and $\phi\in\Phi_4$ be given. We must show
that $(P_1,\ldots,P_N)\in V_\phi$. By the definition of $\Phi_4$ there
exist an index $i_0\in\{1,\ldots,N\}$
and distinct elements $\gamma_1,\ldots,\gamma_4$ of $\Gamma$ such that
$\phi=\{(i_0,\gamma_j):1\le j\le4\}$. 

Since $(P_1,\ldots,P_N)\in \calT^N$, we have $P_{i_0}\in \calT$, and hence $\realcalX_{\Gamma\cdot P_{i_0}}$
is weakly simple. The discrete torsion-free group $\Gamma$ acts freely
on $\HH^3$, and since the elements $\gamma_1,\ldots,\gamma_4$ of
$\Gamma$ are distinct, the points $\gamma_1\cdot
P_{i_0},\ldots,\gamma_4\cdot P_{i_0}$ of $\HH^3$ are distinct; hence 
$
X_{\gamma_1\cdot
  P_{i_0}}^{\Gamma\cdot P_{i_0}},\ldots,
X_{\gamma_4\cdot
  P_{i_0}}^{\Gamma\cdot P_{i_0}}$
are distinct Voronoi regions for 
$\Gamma\cdot P_{i_0}$. By weak simplicity, these four Voronoi regions
cannot share a one-dimensional face, and therefore 
$
X_{\gamma_1\cdot
  P_{i_0}}^{\Gamma\cdot P_{i_0}}\cap\cdots\cap
X_{\gamma_4\cdot
  P_{i_0}}^{\Gamma\cdot P_{i_0}}$ is at most a single point.

Set $\scrS=\scrS(P_1,\ldots,P_N)$.
For any $j\in\{1,\ldots,4\}$, and any point $x\in 
X_{\gamma_j\cdot
  P_{i_0}}^{\scrS}$,
the definition of a Voronoi
region gives $\dist(x,P_{i_0})=\min_{z\in
  \scrS}\dist(x,z)$. Since $P_{i_0}\in\Gamma\cdot
P_{i_0}\subset \scrS$, we have a fortiori that $\dist(x,P_{i_0})=\min_{z\in
  \Gamma\cdot P_{i_0}}\dist(x,z)$, which gives $x\in
X_{\gamma_j\cdot
  P_{i_0}}^{\Gamma\cdot P_{i_0}}$. This shows that
$X_{\gamma_j\cdot
  P_{i_0}}^{\scrS}
\subset
X_{\gamma_j\cdot
  P_{i_0}}^{\Gamma\cdot P_{i_0}}
$ for $i=1,\ldots,4$. Hence
$$
X_{\gamma_1\cdot
  P_{i_0}}^{\scrS}\cap\cdots\cap
X_{\gamma_4\cdot
  P_{i_0}}^{
\scrS
}
\subset
X_{\gamma_1\cdot
  P_{i_0}}^{\Gamma\cdot P_{i_0}}\cap\cdots\cap
X_{\gamma_4\cdot
  P_{i_0}}^{\Gamma\cdot P_{i_0}}.
$$
Since we have seen that
$X_{\gamma_1\cdot
  P_{i_0}}^{\Gamma\cdot P_{i_0}}\cap\cdots\cap
X_{\gamma_4\cdot
  P_{i_0}}^{\Gamma\cdot P_{i_0}}
$
is at most a single point, it follows that
$X_{\gamma_1\cdot
  P_{i_0}}^{\scrS}\cap\cdots\cap
X_{\gamma_4\cdot
  P_{i_0}}^{
\scrS}$ is at most a single point; by definition this
means that $(P_1,\ldots,P_N)\in V_\phi$, as required. This completes
the proof of \ref{bridget}.

Next, for each $\phi\in\Phi$, we denote by $J_\phi$ the set of all points
$(P_1,\ldots,P_N)\in\tU$ for which there exists a circle or
horocycle $C\subset\HH^3$ such that $\gamma\cdot P_i\in C$ for every
$(i,\gamma)\in\phi$. We claim:
\Claim\label{it's closed}
For every $\phi\in\Phi$, the set $J_\phi$ is closed in the subspace
topology of $\tU\subset(\HH^3)^N$.
\EndClaim

To prove \ref{it's closed}, we consider a sequence
$((P_1^{(n)},\ldots, P_N^{(n)}))_{n\ge1}$ of points in $J_\phi$
converging to a point $(P_1^{(\infty)},\ldots,
P_N^{(\infty)})\in\tU$. We must show that $(P_1^{(\infty)},\ldots,
P_N^{(\infty)})\in J_\phi$. For each $n\ge1$ there is a circle or horocycle
$C_n$
containing all the points of the form $\gamma\cdot P_i^{(n)}$ for 
$(i,\gamma)\in\phi$. We use the Poincar\'e model, identifying $\HH^3$
with the unit ball in $\EE^n$. Then each $C_n$ is the intersection of
$\HH^3$ with a Euclidean circle $C_n'$ in the closed unit ball
$\overline{\HH^3}\subset\EE^n$. After passing to a subsequence we may
assume that $(C_n')$ converges, in
the Hausdorff distance for compact subsets in $\EE^3$, to a set
$C_\infty'\subset\overline{\HH^3}$ which is either a Euclidean circle
or a one-point set. If we set
$C_\infty=C'_\infty\cap\HH^3$, then $\gamma\cdot P_i\in C_\infty$ for every
$(i,\gamma)\in\phi$. In particular $C_\infty\ne\emptyset$. This means
that the set $C'_\infty$ is not contained in the unit sphere
$\partial\overline{\HH^3}$, and is therefore either a one-point subset
of
$\HH^3$ or a circle in $\overline{\HH^3}$ intersecting
$\partial\overline{\HH^3}$ in at most one point. But if $C'_\infty$
were a one-point subset $\{A\}$ of $\HH^3$ we would have $\gamma\cdot
P_i^{(n)}=A$ for each 
$(i,\gamma)\in\phi$; this is impossible, since
$(P_1,\ldots,P_N)\in\tU$ and $\Gamma$ acts freely on $\HH^3$.
Hence $C'_\infty$ is a circle in $\overline{\HH^3}$ intersecting
$\partial\overline{\HH^3}$ in at most one point. This implies that
$C_\infty$ is either a circle or a horocycle, and since 
$\gamma\cdot P_i\in C_\infty$ for every
$(i,\gamma)\in\phi$, it follows that $(P_1^{(\infty)},\ldots,
P_N^{(\infty)})\in J_\phi$. Thus \ref{it's closed} is proved.

We set $W_\phi=\tU-J_\phi$ for every $\phi\in\Phi$. Since $\tU$ is
open in $(\HH^3)^N$, 
  \ref{it's closed} immediately implies:
\Claim\label{it's open}
For every $\phi\in\Phi$, the set $W_\phi$ is open in $(\HH^3)^N$.
\EndClaim

Next, we claim:
\Claim\label{so what}
For every $\phi\in\Phi$ we have $W_\phi\subset V_\phi$. 
\EndClaim

To prove \ref{so what}, let $(i_1,\gamma_1),\ldots,(i_4,\gamma_4)$
denote the distinct elements of $\phi$. Suppose that $(P_1,\ldots,P_N)$ is a point of
$(\HH^3)^N$ which does not belong to $V_\phi$. Set
$\scrS=\scrS(P_1,\ldots,P_N)$. By definition, to say that
$(P_1,\ldots,P_N)\notin V_\phi$ 
means that the set
$$E_\phi^{(P_1,\ldots,P_N)}=\bigcap_{j=1}^4X_{\gamma_j\cdot
  P_{i_j}}^{\scrS}$$ contains at least two distinct points 
$Q_1$ and $Q_2$. According to the definition of the Voronoi regions $X_{\gamma_j\cdot
  P_{i_j}}^{\scrS}$, this implies that for $k=1,2$ we have
$$\dist(Q_k,\gamma_1\cdot P_{i_1})=\dist(Q_k,\gamma_2\cdot P_{i_2})=\dist(Q_k,\gamma_3\cdot P_{i_3})=\dist(Q_k,\gamma_4\cdot P_{i_4}).
$$
Hence for $k=1,2$ there is a sphere $S_k$ which is centered at $Q_k$
and contains the points $\gamma_1\cdot P_{i_1},\ldots,\gamma_4\cdot P_{i_4}$. In
particular we have $S_1\cap S_2\ne\emptyset$, and since $S_1$ and
$S_2$ have distinct centers, their intersection is a circle in
$\HH^3$. This circle contains the points $\gamma_1\cdot
P_{i_1},\ldots,\gamma_4\cdot P_{i_4}$, which by definition implies
that $(P_1,\ldots,P_N)\in J_\phi$, and hence that
$(P_1,\ldots,P_N)\notin W_\phi$. This completes the proof of \ref{so what}.

Now we claim:
\Claim\label{it's dense}
For every $\phi\in\Phi_1\cup\Phi_2$, the set $W_\phi$ is dense in $(\HH^3)^N$.
\EndClaim

We have observed that $\tU$ is dense in $(\HH^3)^N$. Hence, in order
to prove \ref{it's dense} it suffices to show that $W_\phi$ is dense
in $\tU$. We will consider an arbitrary point
$(P_1,\ldots,P_N)\in(\HH^3)^N$, and show that $(P_1,\ldots,P_N)$ belongs to the
closure of $W_\phi$ in the subspace topology of $\tU$.

There are two overlapping cases: we have
$\phi\in\Phi_1$ or $\phi\in\Phi_2$. First suppose that
$\phi\in\Phi_1$. Then we may fix an index $i_0\in\{1,2,3,4\}$ such
that $\mu(i_0.\phi)=1$. We may then write the elements of $\phi$ in the
form $(i_0,\gamma_0)$, $(i_1,\gamma_1)$, $(i_2,\gamma_2)$ and
$(i_3,\gamma_3)$, where the indices $i_1$, $i_2$ and $i_3$ are
distinct from $i_0$ (but not necessarily distinct from one
another). 

We may assume that  $(P_1,\ldots,P_N)\in J_\phi=\tU-W_\phi$, as otherwise there is nothing
to prove. Thus there is a
circle or horocycle $C\subset\HH^3$ such that $\gamma_j\cdot P_{i_j}\in C$ for every
$(i,\gamma)\in\phi$.

We can find a point $P_{i_0}'\in\HH^3$, arbitrarily close to
$P_{i_0}$, such that $\gamma_{0}\cdot P_{i_0}'\notin C$. For
$i\in\{1,\ldots,N\}- \{i_0\}$ we set $P_i'=P_i$. 
By taking $P_{i_0}'$ sufficiently close to $P_{i_0}$ we can
guarantee that $(P_1',\ldots,P_N')$ is as close to $(P_1,\ldots,P_N)$
as desired; in particular we can guarantee that
$(P_1',\ldots,P_N')\in\tU$. We shall complete the proof in this case
by showing that $(P_1',\ldots,P_N')\in
W_\phi$. 

For $j=1,2,3$, since
$i_j\ne i_0$, we have $\gamma_{j}\cdot P_{i_j}'=\gamma_{j}\cdot P_{i_j}\in
C$. Furthermore, since we have $(P_1,\ldots,P_N)\in J_\phi\subset\tU$,
the points $\gamma_{1}\cdot P_{i_1}$, $\gamma_{2}\cdot P_{i_2}$ and $\gamma_{3}\cdot P_{i_3}$ are distinct; hence $C$
is the unique circle or horocycle containing $\gamma_{1}\cdot P_{i_1}$, $\gamma_{2}\cdot P_{i_2}$ and $\gamma_{3}\cdot P_{i_3}$. Since $\gamma_{0}\cdot P_{i_0}'\notin C$, it follows
that no circle or horocycle contains all the points
$\gamma_{j}\cdot P_{i_j}'$ for $j=0,1,2,3$. This means that
$(P_1',\ldots,P_N')\notin J_\phi$; since $(P_1',\ldots,P_N')\in\tU$,
we have  $(P_1',\ldots,P_N')\in
W_\phi$. 
This proves the assertion in this case.

Now suppose that $\phi\in\Phi_2$. 
Then we may fix an index $i_0\in\{1,2,3,4\}$ such
that $\mu(i_0.\phi)=2$. We may then write the elements of $\phi$ in the
form $(i_0,\gamma_0)$, $(i_0,\gamma_1)$, $(i_2,\gamma_2)$ and
$(i_3,\gamma_3)$, where the indices $i_2$ and  $i_3$ are
distinct from $i_0$ (but not necessarily distinct from each
other). 

We distinguish two subcases, depending on whether there is a horocycle
containing all the points $\gamma_{0}\cdot P_{i_0}$,
$\gamma_{1}\cdot P_{i_0}$, $\gamma_{2}\cdot P_{i_2}$ and
$\gamma_{3}\cdot P_{i_3}$.

First consider the subcase where there is no such horocycle. 
Again we may assume that  $(P_1,\ldots,P_N)\in J_\phi=\tU-W_\phi$;
then by the definition of $J_\phi$, there is a circle $C$
containing all the points $\gamma_{0}\cdot P_{i_0}$,
$\gamma_{1}\cdot P_{i_0}$, $\gamma_{2}\cdot P_{i_2}$ and
$\gamma_{3}\cdot P_{i_3}$.

Since $\phi$ is by definition a four-element set, the pairs $(i_0,\gamma_0)$, $(i_0,\gamma_1)$, $(i_2,\gamma_2)$ and
$(i_3,\gamma_3)$ are distinct. In particular we have
$\gamma_1\ne\gamma_0$. Since $\Gamma$ is discrete and torsion-free, it
follows that $\delta:=\gamma_0^{-1}\gamma_1$ is
loxodromic or parabolic, and therefore  has no bounded orbits in
$\HH^3$. Since the circle $\gamma_1^{-1}\cdot C$ is compact, $\gamma_1^{-1}\cdot
C$ cannot be
$\delta$-invariant. This implies that the circles $\gamma_0^{-1}\cdot C$ and
$\gamma_1^{-1}\cdot C$ are distinct. Hence $\gamma_0^{-1}\cdot C\cap
\gamma_1^{-1}\cdot C$ consists of at most two points.
Since $\gamma_{1}\cdot P_{i_0}\in C$, it now follows that we can find a point $P_{i_0}'\in\HH^3$, arbitrarily close to
$P_{i_0}$, such that $\gamma_{1}\cdot P_{i_0}'\in C$ but
$\gamma_{0}\cdot P_{i_0}'\notin C$. 
As in the preceding subcase, we set $P_i'=P_i$ for $i\ne i_0$, and
observe that, by taking $P_{i_0}'$ sufficiently close to $P_{i_0}$, we can
guarantee that $(P_1',\ldots,P_N')$ is as close to $(P_1,\ldots,P_N)$
as desired, and in particular that
$(P_1',\ldots,P_N')\in\tU$.

For $j=2,3$, since
$i_j\ne i_0$, we have $\gamma_{j}\cdot P_{i_j}'=\gamma_{j}\cdot P_{i_j}\in
C$. 
Furthermore, since we have $(P_1',\ldots,P_N')\in \tU$,
the points $\gamma_{1}\cdot P_{i_0}'$, $\gamma_{2}\cdot P_{i_2}'$ and $\gamma_{3}\cdot P_{i_3}'$ are distinct; hence $C$
is the unique circle containing $\gamma_{1}\cdot P_{i_0}$, $\gamma_{2}\cdot P_{i_2}$ and $\gamma_{3}\cdot P_{i_3}$, and there is no horocycle
containing them. Since $\gamma_{0}\cdot P_{i_0}'\notin C$, it follows
that no circle or horocycle contains all the points
$\gamma_{0}\cdot P_{i_0}'$,
$\gamma_{1}\cdot P_{i_0}$, $\gamma_{2}\cdot P_{i_2}$ and $\gamma_{3}\cdot P_{i_3}$. This means that
$(P_1',\ldots,P_N')\notin J_\phi$, 
Since $(P_1',\ldots,P_N')\in\tU$,
we have  $(P_1',\ldots,P_N')\in
W_\phi$, and the assertion is proved in this subcase.

There remains the subcase in which there is a horocycle $C$
containing  $\gamma_{0}\cdot P_{i_0}$,
$\gamma_{1}\cdot P_{i_0}$, $\gamma_{2}\cdot P_{i_2}$ and
$\gamma_{3}\cdot P_{i_3}$. We can find a point $P_0'$ arbitrarily
close to $P_0$ but not lying on the horocycle $C$. Again, we set $P_i'=P_i$ for $i\ne i_0$, and
observe that, by taking $P_{i_0}'$ sufficiently close to $P_{i_0}$, we can
guarantee that $(P_1',\ldots,P_N')$ is as close to $(P_1,\ldots,P_N)$
as desired, and in particular that
$(P_1',\ldots,P_N')\in\tU$. 

Since $(P_1,\ldots,P_N)\in\tU$. and since $\Gamma$ acts freely on
$\HH^3$, the points $\gamma_{2}\cdot P_{i_2}=\gamma_{2}\cdot P_{i_2}'$ and
$\gamma_{3}\cdot P_{i_3}=\gamma_{3}\cdot P_{i_3}'$ are distinct. Hence
$C$ is the unique horocycle containing these two points. Since
$P_0'\notin C$, there is no horocycle containing  the points $\gamma_{0}\cdot P_{i_0}'$,
$\gamma_{1}\cdot P_{i_0}'$, $\gamma_{2}\cdot P_{i_2}'$ and
$\gamma_{3}\cdot P_{i_3}'$. It now follows from the preceding subcase
that $(P_1',\ldots,P_N')$ lies in the closure of $W_\phi$. This
completes the proof of \ref{it's dense}.

We now complete the proof that $Y$ contains a dense $G_\delta$ set. It
follows from \ref{dragnet} that 
$$Y=\tU\cap\bigcap_{\phi\in\Phi_1}V_\phi\cap\bigcap_{\phi\in\Phi_2}V_\phi\cap\bigcap_{\phi\in\Phi_4}V_\phi.$$
We have observed that $\tU$ is open and dense in $(\HH^3)^N$. The set
$\bigcap_{\phi\in\Phi_4}V_\phi$ contains a dense $G_\delta$ set according
to \ref{bridget}. 

Now suppose that $m\in\{1,2\}$ is given. The set $\bigcap_{\phi\in\Phi_m}V_\phi$
contains 
$\bigcap_{\phi\in\Phi_m}W_\phi$ by \ref{so what}. The index set
$\Phi_m$ is countable, and each of the sets $V_\phi$ for
$\phi\in\Phi_m$ is open in $(\HH^3)^N$ according to \ref{it's open}, and
is dense in $(\HH^3)^N$ according to \ref{it's dense}. Hence
$\bigcap_{\phi\in\Phi_m}W_\phi$ is a dense $G_\delta$ set. The proof that
$Y$ contains a dense $G_\delta$ set is now complete.
\EndProof

We will need the following elementary fact, Lemma \ref{why proper},
about Voronoi complexes. Recall that a map between locally compact
spaces is said to be {\it proper} if the preimage of every compact set
is compact.

\Lemma\label{why proper}
Let $n\ge2$ be an integer, let $\Gamma$ be a discrete, torsion-free
subgroup of $\HH^n$, let $\tS\subset\HH^n$ be a locally finite,
$\Gamma$-invariant set, let $M$ denote the quotient manifold
$\HH^n/\Gamma$, and let $q:\HH^n\to M$ denote the quotient map. Then
for every polyhedron $X$ in the Voronoi complex $\realcalX_\tS$, the
map $q|X:X\to M$ is proper.
\EndLemma

\Proof
By the definition of $\realcalX_\tS$, the polyhedron $X$ is a face of some
Voronoi region for $\tS$. Hence it suffices to prove the lemma in the
case where $X$ is itself a Voronoi region.

Let $(P_i)_{i\ge1}$ be a sequence of points in $X$ such that the
sequence $(q(P_i))_{i\ge1}$ converges to a point $p\in M$. We must
show that $(P_i)$ has a convergent subsequence. Fix a neighborhood $V$
of $p$ in $M$ which is evenly covered
by the covering map $q$, and fix an open neighborhood $U$ of $p$ whose
closure is contained in $V$ and is compact. We may write $q^{-1}(U)$ as a disjoint union
$\bigcup_{\gamma\in\Gamma}\gamma\cdot\tU$, where $\tU\subset\HH^n$ is
open and has compact closure in $\HH^n$, and $q$ maps $\gamma\cdot \tU$ homeomorphically onto $U$ for every
$\gamma\in\Gamma$.

After passing to a subsequence we may assume
that for every $i$ we have $q(P_i)\in U$ and hence
$P_i\in\gamma_i\cdot\tU$ for some $\gamma_i\in\Gamma$. For each $i$ we
therefore have $X\cap\gamma_i\cdot\tU\ne\emptyset$, and therefore
$\tU\cap\gamma_i^{-1}\cdot X\ne\emptyset$. Since each of the sets
$\gamma_i^{-1}\cdot X$ belongs to the family
$\realcalX_\tS$, which is locally finite by the definition of a
polyhedral complex,
and $\tU$ has compact closure in $\HH^n$, it  now follows that the set
$\{\gamma_i^{-1}\cdot X:i\ge1\}$ is finite. Hence, after passing to a
subsequence, we may assume that the sequence $(\gamma_i^{-1}\cdot
X)_{i\ge1}$ is constant. Thus for any two indices $i$ and $i'$ we have 
$\gamma_{i'}\gamma_i^{-1}\cdot
X=X$. Since the definition of a Voronoi region implies that $X$ contains a unique point of the
$\Gamma$-invariant set $\tS$, 
and since $\Gamma$ acts freely on $\HH^n$, it follows that
$\gamma_{i'}\gamma_i^{-1}$ is the identity. This shows that the
sequence $(\gamma_i)_{i\ge1}$ is constant. We therefore have
$P_i\in\gamma_1\cdot\tU$ for every $i$. Since $(q(P_i))$ converges to
the point
$p\in U$, and $q$ maps $\gamma_1\cdot \tU$ homeomorphically onto $U$, we
deduce that $P_i$ converges in $\gamma_1\cdot\tU$.
\EndProof

\section{Good sets}

\Number\label{gpc review}
In \cite[Definition
2.8]{ratioI}, the notion of  a ``\gpc''  is formally  defined. 
It is a
certain kind of ordered quadruple
$K=(\redcalM,\cala,(\maybecalD_H)_{H\in\cala}, (\Phi_H)_{H\in\cala})$, 
where $\calM$
is a space
and $\cala$ is a locally finite collection of subsets of $\calM$
called cells, which set-theoretically partition $\calM$.
For each  $H\in\cala$,
the
object denoted by
$\calD_H$ is a 
polyhedron in hyperbolic space of some dimension, and
$\Phi_H:\calD_H\to\calM$ is a continuous map
which maps the interior of $\calD_H$ homeomorphically onto $H$; in
particular, $H$ is topologically an open ball of the same dimension as
$\calD_H$. As $H$ varies over $\cala$, the polyhedra $\calD_H$ and the
maps $\Phi_H$ are
subject to a certain natural compatibility condition, which is given
explicitly in \cite{ratioI}. 
(In the special case
where the $\calD_H$ are all compact, the partition $\cala$ in
particular defines a CW complex.)

Now let $n\ge2$ is an integer, and let  $S$ be a non-empty finite subset of
a  hyperbolic $n$-manifold $M$. It is shown in Subsections 2.10 and 2.11 of \cite{ratioI}, and will
briefly reviewed here, how any finite subset $S$ of $M$ gives rise to
a \gpc\ with underlying space $M$. 
We
write $M=\HH^n/\Gamma$, where $\Gamma\le\isom(\HH^n)$ is discrete,
and torsion-free, and let $q:\HH^n\to
M
$ denote the
quotient  map. Then $\tS:=q^{-1}(S)$ is a non-empty, locally finite
subset of $\HH^n$, and hence by the material reviewed in \ref{voronoi
  review} above, 
$\realcalX_{\tS}$ is a well-defined \willbepolytopal\
complex with underlying space $\HH^n$. It is 
$\Gamma$-invariant in the sense that $\gamma\cdot
X\in\realcalX_\tS$ whenever $\gamma\in\Gamma$ and
$X\in\realcalX_\tS$. We define a partition 
$\cala_S$
of $M$ to consist of all sets of the form $q(\inter X)$ where $X$ is a
polyhedron in $\realcalX_{\tS}$. 
For each $H\in\cala_S$
we choose
an
element $X$ of $\realcalX_{\tS}$ with 
$H=q(\inter X)$, 
and we set $\calD_{H,S}=X$ and $\Phi_{H,S}=q|X$. 
In Subsections 2.10 and 2.11 of \cite
{ratioI} it is observed that
the set $S$ determines the
partition $\cala_S$; that it determines each $\calD_{H,S}$ up to
isometry; that it determines each
$\Phi_{H,S}$ up to precomposition with an
isometry between convex polyhedra; and that 
$K_S:=(M,\cala_S,(\maybecalD_{H,S})_{H\in\cala_S}, (\Phi_{H,S})_{H\in\cala_S})$, 
is a \gpc.

Let $X$ be any $3$-dimensional polyhedron in
$\realcalX_\tS$. Then it follows from the constructions reviewed above that $q$ maps $\inter X$ homeomorphically onto a
$3$-cell $H$ of $K_S$, and that there is a unique isometry $J_{X,S}:
X\to\calD_{H,S}$ of convex hyperbolic polyhedra such that 
$\Phi_{H,S}\circ J_{X,S}=q|X$.
\EndNumber

\Number\label{thick and thin}
As in \cite{ratioI}, if $M$ is a hyperbolic $3$-manifold and
$\epsilon$ is a positive real number, we will
denote by $\Mthin(\epsilon)$ the set of points $p\in M$ which are {\it
  $\epsilon$-thin} in the sense that there is a homotopically
non-trivial closed path based at $p$ whose length is less than
$\epsilon$; and we set $\Mthick(\epsilon)=M-\Mthin(\epsilon)$. 

At a few
points in this paper, it will be convenient to use the following
notation which was used in \cite {kfree-volume}. If $p$ is a point of a
non-simply connected hyperbolic $3$-manifold $M$, we denote by
$\shortone(p)$ the length of the shortest homotopically non-trivial
closed path in $M$ based at $p$. Thus $\shortone$ is a positive-valued
function on $M$, and for any $\epsilon>0$ we have
$\Mthin(\epsilon)=\{p\in M:\shortone(p)<\epsilon\}$ and
$\Mthick(\epsilon)=\{p\in M:\shortone(p)\ge\epsilon\}$.

A {\it Margulis number,} for an orientable hyperbolic $3$-manifold $M$
is defined to be a positive number $\tryepsilon$ such that, for every point
$p\in M$ and for any two closed paths $\alpha$ and $\beta$  based at $p$
such that the elements $[\alpha]$ and $[\beta]$ of $\pi_1(M,p)$ do not
commute, we have $\max(\length\alpha,\length\beta)\ge\tryepsilon$. As
in \cite{ratioI}, if the
strict inequality $\max(\length\alpha,\length\beta)>\tryepsilon$ holds for all
$\alpha$ and $\beta$ such that $[\alpha]$ and $[\beta]$ do not
commute, we will say that $\tryepsilon$ is a {\it strict Margulis number} for
$M$.

According to \cite[Proposition 4.5]{ratioI}, if $\epsilon$ is a strict
Margulis number for $M$, then $\Mthick(\epsilon)$ is a smooth $3$-manifold-with-boundary. 
\EndNumber

\Definition\label{good-def}
Let $M$ be an orientable hyperbolic $3$-manifold, and let $\epsilon$
be a strict Margulis number for $M$.
Let $S\subset M$ be a non-empty finite set. Let us
write $M=\HH^3/\Gamma$, where $\Gamma\le\isomplus(\HH^3)$ is discrete
and torsion-free, and set $\tS=q^{-1}(S)$, where $q:\HH^3\to M$
denotes the quotient map. We will say that $S$ is  {\it $\epsilon$-\good} 
if 
(1) the polyhedral complex $\realcalX_\tS$ is weakly simple (see \ref{weakly simple def.}), and (2) every
cell of $K_S$ which meets $\Mthick(\epsilon)$ meets the interior of
the $3$-manifold-with-boundary
$\Mthick(\epsilon)$ (see \ref{thick and thin}).
\EndDefinition

\Number\label{nets review}
As in \cite{ratioI}, we define an {\it $\epsilon$-\net} for a metric
space $X$, where $\epsilon$ is a positive real
number, to be a subset $S$ of
$X$ such that for any two distinct elements $x,y$ of $S$ we have
$\dist(x,y)\ge\epsilon$. We also recall from \cite[Subsection
4.6]{ratioI} that
an {\it $\maybeepsilon$-thick $\maybeepsilon$-\net} for a hyperbolic manifold
$M$ is defined to be an $\maybeepsilon$-\net\ for $M$ which is contained in
$\Mthick(\maybeepsilon)$; and that by a {\it maximal $\maybeepsilon$-thick $\maybeepsilon$-\net} for
$M$ one simply means an $\maybeepsilon$-thick $\maybeepsilon$-\net\ for
$M$ which is not properly contained in any other $\maybeepsilon$-thick $\maybeepsilon$-\net\ for
$M$.
\EndNumber

\Lemma\label{when is it}
Let $M$ be a finite-volume orientable hyperbolic $3$-manifold, and let $\epsilon_0$
be a Margulis number for $M$. Then 
for every positive number
$\epsilon_1<\epsilon_0$, there is a number $\epsilon$ 
such that
\begin{itemize}
\item
$\epsilon_1<\epsilon<\epsilon_0$ (so that $\epsilon$ is a strict
Margulis number for $M$), and
\item
$M$ contains a maximal
$\epsilon$-thick $\epsilon$-\net\ which is $\epsilon$-\good.
\end{itemize}
\EndLemma

\Proof
We write $M=\HH^3/\Gamma$, where $\Gamma\le\isomplus(\HH^3)$ is
discrete, and we let $q:\HH^3\to M$ denote the quotient map.

Since $M$ has finite volume, it follows from 
\cite[Proposition 4.5]{ratioI}
that for every $\alpha>0$, the space $\Mthick(\alpha)$
is compact.
Hence if we regard $\Mthick(\alpha)$ as a
metric space with the 
``extrinsic'' distance function obtained by restricting the hyperbolic
distance function on $M$,
there is an upper
bound on the cardinality of any $\alpha$-\net\ in $\Mthick(\alpha)$. 
For every $\alpha>0$ we      define a non-negative integer $\nu(\alpha)$ to be
the largest cardinality of any $\alpha$-thick $\alpha$-\net\ in
$M$. 
Then $\nu$ is a (weakly) monotone decreasing integer-valued
function on $(0,\infty)$. Hence there is a number $\epsilon_2$ with
$\epsilon_1\le\epsilon_2<\epsilon_0$ such that $\nu$ is constant on
the half-open interval $[\epsilon_2,\epsilon_0)$. Let $N$ denote the
constant value of $\nu$ on this interval.

Let us choose numbers $\epsilon_3$ and $\epsilon_4$ with
$\epsilon_2<\epsilon_4<\epsilon_3<\epsilon_0$. We have
$\nu(\epsilon_3)=N$, and hence there is an $\epsilon_3$-thick
$\epsilon_3$-\net\ $S_0\subset M$ with $\#(S_0)=N$. We write
$S_0=\{p_1^0,\ldots,p_N^0\}$, and we choose points
$P_1^0,\ldots,P_N^0\in\HH^3$ such that $q(P_i^0)=p_i^0$ for
$i=1,\ldots,N$, According to Proposition \ref{when would it be},
there exist points $P_1,\ldots,P_N\in\HH^3$, arbitrarily close to
$P_1^0,\ldots,P_N^0$ respectively, such that the following conditions hold.
\Claim\label{twiddle}
\textnormal{The orbits $\Gamma\cdot P_1,\ldots,\Gamma\cdot P_N$ are distinct.}
\EndClaim
 \Claim\label{twaddle}
\textnormal{If we denote by $\scrS$ the locally finite family $\Gamma\cdot
  P_1\cup\cdots\cup\Gamma\cdot P_N$, the polyhedral complex
  $\realcalX_\scrS$ is weakly simple.}
\EndClaim
We set $p_i=q(P_i)$ for $i=1,\ldots,N$. Then  \ref{twiddle} says
that $p_1,\ldots,p_N$ are distinct, so that $S:=\{p_1,\ldots,p_N\}$
has cardinality $N$. Note that the set denoted $\scrS$ in \ref{twaddle} is $q^{-1}(S)$.

Since $S_0$ is an $\epsilon_3$-thick
$\epsilon_3$-\net, we have $\shortone(p_i^0)\ge\epsilon_3$ for
$i=1,\ldots,N$
(in the notation reviewed in \ref{thick and thin}),
and $\dist(p_i^0,p_j^0)\ge\epsilon_3$ whenever $i\ne
j$.
By choosing the $P_i$ sufficiently close to the $P_i^0$ we can
guarantee that $p_i$ is as close as desired to $p_i^0$ for
$i=1,\ldots,N$; since the distance function on the metric space $M$ is
continuous, and since $\shortone$ is continuous by 
\cite[Lemma 4.5]{kfree-volume},
we may suppose the $P_i$ to be chosen so that 
$\shortone(p_i)>\epsilon_4$ for
$i=1,\ldots,N$, and $\dist(p_i,p_j)>\epsilon_4$ whenever $i\ne
j$. Thus  $S$ is an $\epsilon_4$-thick
$\epsilon_4$-\net.

For each $C\in\redasubs$ 
(see \ref{gpc review}),
we set
$m_C=\sup_{p\in C}\shortone(p)\in\RR\cup\{+\infty\}$. The family
$\redasubs$ is locally finite; this is by definition included in the
statement that $K_S$ is a \gpc, and is directly verified in
Subsections 2.10 and 2.11 of \cite{ratioI}. In particular
$\redasubs$ is countable, and hence the set
$E:=\{m_C:C\in\redasubs\}\subset\RR\cup\{+\infty\}$ is countable. We
may therefore choose a real number
$\epsilon$ with $\epsilon_2<\epsilon<\epsilon_4$ such that
$\epsilon\notin E$. In particular we have
$\epsilon_1<\epsilon<\epsilon_0$ 
(so that $\epsilon$ is a strict
Margulis number for $M$).

Since $\epsilon<\epsilon_4$, the 
$\epsilon_4$-thick
$\epsilon_4$-\net\ $S$ is an $\epsilon$-thick
$\epsilon$-\net. Since $\epsilon_2<\epsilon<\epsilon_0$, we have
$\nu(\epsilon)=N=\#(S)$, which implies that the $\epsilon$-thick
$\epsilon$-\net\ $S$ is maximal. It remains to show that $S$ is $\epsilon$-good.

According to \ref{twaddle}, $S$ satisfies Condition (1) of 
Definition \ref{good-def}. To check Condition (2) of that definition,
suppose that $C$ is a cell of $K_S$
such that
$C\cap\Mthick(\epsilon)\ne\emptyset$. Thus $\shortone$ takes a value
greater than or equal to $\epsilon$ at some point of $C$, and in
particular $m_C=\sup_{p\in C}\shortone(p)\ge\epsilon$. But since
$\epsilon\notin E$, we have $m_C\ne\epsilon$; hence $m_C>\epsilon$. This
means that $\shortone$ takes a value strictly
greater than $\epsilon$ at some point 
$p_0\in C$. Since the function $\shortone$ is continuous on $M$ by
\cite[Lemma 4.5]{kfree-volume}, the point $p_0$ has a neighborhood in $M$ on which $\shortone$
is bounded below by $\epsilon$; such a neighborhood is contained in 
$\Mthick(\epsilon)$ according to \ref{thick and thin}. This shows that
$p_0\in\inter(\Mthick(\epsilon))$, so that
$C\cap\inter(\Mthick(\epsilon))\ne\emptyset$. This gives
Condition (2).
\EndProof

\Lemma\label{before and what if it is}
Let $M$ be  a 
finite-volume 
orientable hyperbolic $3$-manifold, 
let $\epsilon$ be a strict Margulis number for $M$, and let $S\subset
\Mthick(\epsilon)$ be a non-empty finite set which is
$\epsilon$-\good. Let us write $M=\HH^3/\Gamma$, where $\Gamma\le\isomplus(\HH^3)$ is discrete
and torsion-free, with finite covolume, and let $q:\HH^3\to M$ denote the
quotient map. Set $\tS=q^{-1}(S)$, and let $F$ be any polyhedron of
dimension at least $2$, belonging to $\realcalX_\tS$, such that $F\cap q^{-1}(\Mthick(\epsilon))\ne\emptyset$.
Then $F$ has
a proper face
$L_0$  such that $ L_0\cap
q^{-1}(\Mthick(\epsilon))\ne\emptyset$. 
\EndLemma

\Proof
According to Lemma \ref{why proper},
the map $q|F:F\to M$ is proper. 
But since $\vol M<\infty$, the set
$\Mthick(\epsilon)$ is compact by
\cite[Proposition 4.5]{ratioI}.
The properness of $q|F$ then implies that $F\cap
q^{-1}(\Mthick(\epsilon))$ is compact. 
Applying Lemma \ref{cincinnati},
taking $C=F\cap
q^{-1}(\Mthick(\epsilon))\subset F$, 
we deduce that
there is a
compact convex polyhedron $F_0\subset F$, with $\dim F_0=\dim F$,
such that 
$F\cap
q^{-1}(\Mthick(\epsilon))\subset F_0$, 
and
each proper face of $F_0$ is
contained either in a proper face of $F$ or in
$q^{-1}(\Mthin(\epsilon))$.

Assume that the conclusion of the lemma is false, i.e. that every
proper face of $F$ is contained in $q^{-1}(\Mthin(\epsilon))$. Then our
choice of $F_0$ implies that every
proper face of $F_0$ is contained in $q^{-1}(\Mthin(\epsilon))$.
Since the convex polyhedron $F_0$ is compact, the union $\partial F_0$
of the proper faces of $F_0$ is homeomorphic to $S^{d-1}$ 
where $d:=\dim F_0=\dim F\ge2$ (see  \cite[Subsection 2.1.1]{ratioI}).
Hence $\partial F_0$ is
connected, and is therefore contained in a component $Z$ of
$q^{-1}(\Mthin(\epsilon))$. But since $\epsilon$ is a strict Margulis
number for $M$, 
the last sentence of \cite[Proposition 4.5]{ratioI}
asserts that each component of $q^{-1}(\Mthin(\epsilon))$ is 
convex. The compact convex polyhedron $F_0$ is
the convex hull of its boundary, and is therefore contained in $Z$;
hence 
$F_0\subset q^{-1}(\Mthin(\epsilon))$. 
Since 
$F\cap
q^{-1}(\Mthick(\epsilon))\subset F_0$, 
it follows that 
$F\cap q^{-1}(\Mthick(\epsilon))=\emptyset$, 
a contradiction to the
hypothesis.
\EndProof

\RemarkDefinitions\label{one way}
\textnormal{
Let $S$ be a non-empty finite subset of a  hyperbolic $3$-manifold. We
recall from \cite[Definition 2.12]{ratioI} that a {\it \dotsystem}
for $S$ is defined 
to be a subset  $\redcalT $ of $|\atwoS|\subset K_S$ which contains at most one point
of each  $2$-cell of $K_S$.
}

\textnormal{\noindent
Now let $M$ be
  an orientable
  hyperbolic $3$-manifold, and 
let $\epsilon$ be
a strict Margulis
number for $M$.
If $S$ is any non-finite subset of $M$, we may obtain a \dotsystem\ for $S$ by the following construction.
Let $\calc\subset\atwoS$ denote the set of all  $2$-cells of $K_S$
which meet the interior of the $3$-manifold-with-boundary
$\Mthick(\epsilon)$ (see \ref{thick and thin}).
For each $C\in\calc$, select a point $\tau_C\in C\cap\inter(\Mthick(\epsilon))$, and set $\redcalT =\{\tau_C:C\in\calc\}$. That $\redcalT $ is a \dotsystem\ for $S$ is immediate from the definition.
}

\textnormal{\noindent
A \dotsystem\ $\redcalT $ for $S$ that is obtained from the construction
described above will be termed {\it $\epsilon$-\adapted.}
}
\EndRemarkDefinitions

\Number\label{new graph stuff}
By a {\it graph} we mean a CW complex of dimension at most $1$.
The $0$-cells and $1$-cells of a graph will be called its {\it
  vertices} and {\it edges} respectively. 
The underlying space of a graph $\redG$ will be denoted $|\redG|$.
(In the more informal language of \cite{ratioI}, we said that a space
$\scrG$ has the structure of a graph to mean that it is identified
with the underlying space of a graph.)

A {\it parametrization} of
an edge $e$ of a graph $\redG$ is defined to be a path
$\alpha:[0,1]\to|\redG|$ such that $\alpha$ maps $(0,1)$
homeomorphically onto $e$.
(Thus the map $x\mapsto\alpha((x+1)/2)$ from $D^1=[-1,1]$ to $|G|$ is
a characteristic map for the $1$-cell $e$.)

An {\it orientation} of an edge $e$ is an equivalence class of
parametrizations of $e$, where two parametrizations are defined to be
equivalent if one is obtained from the other by precomposition with a monotone
increasing self-homeomorphism of $[0,1]$. An {\it oriented edge}
$\eta$ is defined by an edge $e$ and an orientation of $e$; we call
$e$ the underlying edge of $\eta$. 
By a {\it parametrization of an oriented edge $\eta$} we mean a
parametrization of the underlying edge $e$ of $\eta$ which represents the
orientation of $e$ that defines $\eta$.
If $\eta$ is an oriented edge, we
denote by $\bareta$ the oriented edge which has the same underlying
edge 
as $\eta$,
but has 
the opposite orientation.

The {\it initial vertex} and {\it terminal vertex} of an oriented edge $\eta$,
denoted respectively by $\init_\eta$ and $\term_\eta$, are the vertices
$\alpha(0)$ and $\alpha(1)$, where $\alpha$ is a(n arbitrary)  parametrization of
the oriented edge $\eta$. 
For any oriented
edge $\eta$ we have 
$\init_{\bareta}=\term_\eta$  and $\term_{\bareta}=\init_\eta$.
We call $\eta$ an
oriented {\it loop}, and call its underlying edge a loop, if
$\init_\eta=\term_\eta$. 
\EndNumber

\Number\label{graph review}
Let $S$ be a non-empty finite subset of a  hyperbolic $3$-manifold,
and let $\redcalT $ be a \dotsystem\ for $S$. 
As in \cite[Subsection 2.13]{ratioI},
we will denote by
$\redscrPST$
the set of all
ordered pairs $(H,\wasU)$ such that $H\in\athreeS$ and
$\wasU\in{\Phi}_{H,S}^{-1}(\redcalT )\subset{\Phi}_{H,S}^{-1}(\atwoS)\subset\partial
\maybecalD_{H,S} 
$.
There is a map $\redscrPST\to \redcalT $ defined by
$(H,\wasU)\mapsto\Phi_{H,S}(\wasU)$. According to \cite[Subsection 2.13]{ratioI}, each fiber of this map has cardinality exactly
$2$. 
In this paper we shall denote by $\sigma_{S,\redcalT }$ the involution of $\redscrPST$ that
interchanges the points of each fiber.

It is pointed out in \cite[Subsection 2.13]{ratioI} that for each
$H\in\athreeS$   there is a unique point $\willbep_{H,\willbeS}$ of $S$
lying in $H$, and that there is a unique point $\willbehatp_{H,\willbeS}$  
of $\inter \maybecalD_{H,\willbeS}$ which is mapped to $\willbep_{H,\willbeS}$ by 
$\Phi_{H,S}$. When the set $S$
is understood,
 we will write  $p_{H}$ and  $\hatp_{H}$ 
in place of $p_{H,S}$ and  $\hatp_{H,S}$. As in \cite[Subsection 2.13]{ratioI}, for
each 
$(H,\wasU)\in\scrP_{
S,\redcalT }
$, we will denote by $\ell_{H,\wasU}$ the line segment in the
convex \willbepolytope\ $ \maybecalD_{H,\willbeS}$ joining the points
$\willbehatp_{H,\willbeS}$ and $\wasU$; and
we will denote by
$\scrG^{S,\redcalT }$ the set $S\cup\bigcup_{(H,\wasU)\in
\scrP_{S,\redcalT }
}\Phi_{H,S}(\ell_{H,\wasU})\subset
M$.

It is observed in \cite[Subsection 2.13]{ratioI} (in slightly less formal
language) that there is a naturally
defined graph $G^{S,\redcalT }$ with $|G^{S,\redcalT }|=\scrG^{S,\redcalT }$. The vertex set
of $G^{S,\redcalT }$ is $S$. The edge set of $G^{S,\redcalT }$ is in natural bijective
correspondence with $\redcalT $. If a point $\tau\in \redcalT $ is given, and if the
points of the fiber of $\tau$ under the map 
$(H,\wasU)\mapsto\Phi_{H,S}(\wasU)$ from $\scrP_{S,\redcalT }$ to $\redcalT $ are
denoted by $(H_0,\wasU_0)$ and $(H_1,\wasU_1)$ (so that $\sigma_{S,\redcalT }$
interchanges $(H_0,\wasU_0)$ and $(H_1,\wasU_1)$), then the edge of
$\scrG^{S,\redcalT }$ corresponding to $\tau$ is $\{\tau\}
\cup\Phi_{H_1,S}(\inter\ell_{H_1,\wasU_1})\cup\Phi_{H_2,S}(\inter\ell_{H_2,\wasU_2})$.

We will denote
by $\veccale^{S,\redcalT }$ the set of all oriented edges of the graph
$\redG^{S,\redcalT }$.
\EndNumber

\RemarksNotation\label{before adjacency}
Suppose that $S$ is a non-empty finite subset of a  hyperbolic
$3$-manifold, and that $\redcalT $ is a \dotsystem\ for $S$. We set
$\sigma=\sigma_{S,\redcalT }$. For each
$(H,\wasU)\in\redscrPST$, we will denote by $\lambda_{H,\wasU}$ 
the path which parametrizes the geodesic arc $\ell_{H,\wasU}\subset\calD_{H,S}$ proportionally to arclength
and has initial point $\hatp_H$
and terminal point $\wasU$. We denote by $\xi_{H,\wasU}$ the path 
$\Phi_{H,S}\circ\lambda_{H,\wasU}$ in $M$. If we set $(\Hstar,\Ustar)=\sigma(H,\wasU)$, the path
$\xi_{H,\wasU}\star\overline{\xi_{\Hstar,\Ustar}}$ is then a
parametrization of an edge of $G^{S,\redcalT }$, and in particular it defines
an oriented edge of $G^{S,\redcalT }$, which will be denoted $\eta_{H,\wasU}\in \veccale^{S,\redcalT }$. 
The
assignment $(H,\wasU)\mapsto\eta_{H,\wasU}$ is a bijection from $\redscrPST$ to
$\veccale^{S,\redcalT }$.
We will denote the inverse
assignment, which is a map from $\veccale^{S,\redcalT }$ to
$\redscrPST$, by $\eta\mapsto(\Hinit^\eta,\Uinit^\eta)$. 
Thus for any   $\eta\in\veccale^{S,\redcalT }$ 
we have 
by definition that $(\Hinit^\eta,\Uinit^\eta)\in\redscrPST$; that is,
$\Hinit^\eta$ is a $3$-cell of $K_S$, and $\Uinit^\eta$ is a point of
$\partial\calD_{\Hinit^\eta,S}$ which is mapped by $\Phi_{\Hinit^\eta,S}$ to a point of
$\redcalT $. In particular 
 $\Uinit^\eta$ is an
interior point of a two-dimensional face of $\calD_{\Hinit^\eta,S}$.
This face will be denoted by $\Cinit^\eta$.

If $(H,\wasU)\in\redscrPST$ is given, and we write $(\Hstar,\Ustar)=\sigma(H,\wasU)$, 
we
have $\eta_{\Hstar,\Ustar}=\overline{\eta_{H,\wasU}}$. Hence for any
$\eta\in\veccale^{S,\redcalT }$ we have $(\Hinit^{\bareta},\Uinit^{\bareta})=\sigma(\Hinit^\eta,\Uinit^\eta)
$.

We observe that for every
 $\eta\in\veccale^{S,\redcalT }$ we have $p_{\Hinit^\eta}=\init_\eta$.

For every $\eta\in\veccale^{S,\redcalT }$ we define $\Hterm^\eta=\Hinit^{\bareta}$, $\Uterm^\eta=\Uinit^{\bareta}$, and
$\Cterm^\eta=\Cinit^{\bareta}$. It now follows that 
$(\Hterm^{\eta},\Uterm^{\eta})=\sigma(\Hinit^\eta,\Uinit^\eta)
$; that
 $\Cterm^\eta$ is the
two-dimensional face of $\Hterm^\eta$ whose interior
contains the point $\Uterm^\eta$; and that
$p_{\Hterm^\eta}=\term_\eta$.

\EndRemarksNotation

\Lemma\label{match 'em}
Let $S$ be a non-empty finite subset of a hyperbolic
$3$-manifold $M$. Let $\redcalT $ be a
\dotsystem\ for $S$. Write $M=\HH^3/\Gamma$, where $\Gamma\le\isom(\HH^3)$ is discrete
and torsion-free,  let $q:\HH^3\to
M
$ denote the
quotient  map, and set $\tS=q^{-1}(S)$, Suppose that $\eta\in\veccale^{S,\redcalT }$ is given, that  $\alpha$ is a  parametrization of
the oriented edge $\eta$,
and that
$\talpha:[0,1]\to\HH^3$ is a lift of $\alpha$. For $i=0,1$, set
$P_i=\talpha(i)$ and $X_i=X_{P_i}\in\realcalX_\tS$ (in the notation of \ref{voronoi review}). Then 
we have $q(\inter X_i)=H_i^\eta$ for $i=0,1$.
Furthermore, $X_0\cap X_1$ is a common
two-dimensional face $F$ of $X_0$ and $X_1$; and for $i=0,1$, the isometry $J_{X_i,S}:
X_i\to\calD_{H_i^\eta,S}$ (see \ref{gpc review}) maps the face $F$ of
$X_i$ onto $\wasC_i^\eta$. 
\EndLemma

\Proof
For $i=0,1$ we set $H_i=H_i^\eta$, $u_i=u_i^\eta$, $U_i=U_i^\eta$,
$p_i=p_{H_i}$, $\Phi_i=\Phi_{H_i,S}$, $\calD_i=\calD_{H_i^\eta,S}$, and $J_i=J_{X_i,S}:
X_i\to\calD_i$. According to the discussion in \ref{graph review}, the
elements $(H_0,u_0)$ and $(H_1,u_1)$ of $\scrP_{S,\redcalT }$ are distinct
(and are interchanged by the involution $\sigma_{S,\redcalT }$).

According to \ref{before adjacency} we have
$\init_\eta=p_0$ and $\term_\eta=p_1$. Hence for $i=0,1$ we
have $\alpha(i)=p_i$, and therefore $q(P_i)=p_i$. Now $q(\inter X_i)$  is a  $3$-cell of $K_S$ by
\ref{gpc review}; and since $P_i\in\inter(X_i)$, and $H_i$ is the
unique $3$-cell of $K_S$ containing $p_i=q(P_i)$, we deduce that
$q(\inter X_i)=H_i$.

Using the notation of \ref{before adjacency}, we set
$\lambda_i=\lambda_{H_i,u_i}$ and $\xi_i=\xi_{H_i,u_i}$ for
$i=0,1$. We may assume without loss of generality that $\alpha$ is the
standard parametrization $\xi_0\star\overline{\xi_1}$ of the oriented
edge $\eta$. Since $\talpha$ is a lift of $\alpha$ having initial
point $P_0$ and terminal point $P_1$, there are lifts $\txi_0$ and
$\txi_1$ of $\xi_0$ and $\xi_1$, respectively, such that $\txi_i(0)=P_i$
for $i=0,1$, and such that $\txi_0(1)=\txi_1(1)$; and we have
$\talpha=\txi_0\star\overline{\txi_1}$. But for $i=0,1$ since
$\Phi_i\circ J_i=q|X_i$ by \ref{graph review}, the path
$J_i^{-1}\circ\lambda_i$ in $X_i\subset\HH^3$ is a lift of $\xi_i$. The initial
point of this lift lies in $X_I$ and is mapped by $q$ to $p_i$, and is
therefore equal to $P_i$. 
By the uniqueness of path-lifting in covering spaces, it follows that
$\txi_i=J_i^{-1}\circ\lambda_i$ for $i=0,1$. In particular, 
the paths $J_0^{-1}\circ\lambda_0$ and $J_1^{-1}\circ\lambda_1$ have the same
terminal point, i.e. $J_0^{-1}(u_0)=J_1^{-1}(u_1)$. Since $U_i$ is the unique
two-dimensional face of $H_i$ containing $u_i$, we must have
$J_0^{-1}(U_0)=J_1^{-1}(U_1)=F$, where $F$ is a common two-dimensional face of
$X_0$ and $X_1$. Thus $J_i(F)=U_i$ for $i=0,1$. 

It remains only to show that $X_0\cap X_1=F$. Since the polyhedral complex $\realcalX_\tS$
has underlying space $\HH^3$,
and 
the 
three-dimensional polyhedra $X_0,X_1\in 
\realcalX_\tS$ share the two-dimensional face $F$, we must have either $X_0\cap X_1=F$ or
$X_0=X_1$. Assume that $X_0=X_1$. Then we have
$H_0=q(\inter(X_0))=q(\inter(X_1))=H_1$. The equality $H_0=H_1$  implies
that $J_0=J_1$, and hence that $U_0=J_0(F)=J_1(F)=U_1$; and since
$u_i$ is the unique point of $U_i\cap\Phi_i^{-1}(\redcalT )$ for $i=0,1$, we have
$u_0=u_1$. This is a contradiction, since the elements $(H_0,u_0)$ and $(H_1,u_1)$ of $\scrP_{S,\redcalT }$ are distinct.
\EndProof

\DefinitionRemark\label{adjacency def}
Let $S$ be a non-empty finite subset of a  hyperbolic
$3$-manifold, and let $\redcalT $ be a \dotsystem\ for $S$.

Let $\eta,\eta'\in\veccale^{S,\redcalT }$ be given. If $\init_\eta=\init_{\eta'}$, then 
in the notation of \ref{before
  adjacency} we have $H:=\Hinit^\eta=\Hinit^{\eta'}$, and 
   $\Cinit^\eta$ and $\Cinit^{\eta'}$ are two-dimensional faces of $\maybecalD_{H ,S}$.
We shall say that $\eta$
and $\eta'$ are {\it initially adjacent} if 
   $\Cinit^\eta\cap\Cinit^{\eta'}$ is a one-dimensional face of $\maybecalD_{H ,S}$.

Likewise, if 
$\eta,\eta'\in\veccale^
{S,\redcalT }$ 
and
$\term_\eta=\term_{\eta'}$, then 
$H:=\Hterm^\eta=\Hterm^{\eta'}$, and 
 $\Cterm^\eta$ and $\Cterm^{\eta'}$ are two-dimensional faces of $\maybecalD_{H ,S}$.
We shall say that $\eta$
and $\eta'$ are {\it terminally adjacent} if 
   $\Cterm^\eta\cap\Cterm^{\eta'}$ is a one-dimensional face of $\maybecalD_{H ,S}$.

Note that if $\eta,\eta'\in\veccale^{S,\redcalT }$ satisfy $\init_\eta=\init_{\eta'}$, then
$\term_{\bareta}=\term_{\bareta'}$; and $\bareta$ and $\bareta'$ are
 terminally adjacent if and only if $\eta$ and $\eta'$ are initially adjacent.
\EndDefinitionRemark

\Lemma\label{and what if it is}
Let $M$ be  a 
finite-volume 
orientable hyperbolic $3$-manifold, 
let $\epsilon$ be a strict Margulis number for $M$, and let $S\subset \Mthick(\epsilon)$ be a non-empty finite set which is $\epsilon$-\good.
Let $\redcalT $ be a \dotsystem\ for $S$
which is $\epsilon$-\adapted.
Suppose that $\eta$ is an oriented edge of $\redG^{S,\redcalT }$. Then there
exist oriented edges $\reddenedetanought,\reddenedetaone$ of  $\redG^{S,\redcalT }$ 
such that
\begin{enumerate}
\item $\init_{\reddenedetanought}=\init_\eta$, and
$\reddenedetanought$ is initially adjacent to $\eta$;
\item $\term_{\reddenedetaone}=\term_\eta$, and
$\reddenedetaone$ is terminally adjacent to $\eta$;
\item $\term_{\reddenedetanought}=\init_{\reddenedetaone}$; and
\item if $\alpha$ is a parametrization of the oriented edge
  $\eta$, and if $\alpha_i$ is a parametrization of the oriented edge
  $\eta_i$ for 
$i=0,1$, 
then 
$\alpha_0\star\alpha_1$ 
is fixed-endpoint homotopic to $\alpha$.
\end{enumerate}
\EndLemma

\Proof
We write $M=\HH^3/\Gamma$, where $\Gamma$ is discrete and
torsion-free, and let $q:\HH^3\to M$ denote the quotient map. Set
$\tS=q^{-1}(S)$.

Choose
a parametrization $\alpha$ of the oriented edge $\eta$, and choose a lift $\talpha$ of $\alpha$
to the covering space $q:\HH^3\to M$. 

For $i=0,1$, set
$P_i=\talpha(i)$ and $X_i=X_{P_i}\in\realcalX_\tS$ (see \ref{voronoi review}). 

For $i=0,1$ we set $H_i=H_i^\eta$ and $\calD_i=\calD_{H_i,S}$. We set
$p_i=p_{H_i}=q(P_i)$, so that $\init_\eta=p_0$ and $\term_\eta=p_1$. We also
set $J_i=J_{X_i,S}:
X_i\to\calD_i$; 
$\Phi_i=\Phi_{H_i,S}$; and
$\wasC_i=\wasC_i^\eta$.
According to 
Lemma
\ref{match 'em},  
we have
$q(\inter X_i)= 
H_i$ for $i=0,1$.
Furthermore, $F:=X_0\cap X_1$ is a common
two-dimensional face of $X_0$ and $X_1$, and $J_i(F)=
\wasC_i$ 
for $i=0,1$.

By definition (see \ref{before adjacency}), 
$\wasC_0$ contains the point $\wasU_0^\eta\in
\Phi_0^{-1}(\redcalT )$.
But since the \dotsystem\ $\redcalT $ is $\epsilon$-\adapted, we have 
$\redcalT  \subset \inter(\Mthick(\epsilon)) \subset \Mthick(\epsilon)$ (see \ref{one way}).
It follows that $\wasC_0\cap
\Phi_0^{-1}(\Mthick(\epsilon))\ne\emptyset$, and therefore
$ F\cap  q^{-1}(\Mthick(\epsilon))\ne\emptyset$. 
Hence, 
by Lemma \ref{before and what if it is}, we may fix a proper face
$Y$ of $F$ such that $ Y\cap
q^{-1}(\Mthick(\epsilon))\ne\emptyset$.

Since $\dim F=2$, the proper face $Y$ of $F$ is a (not necessarily
proper) face of a
$1$-dimensional face $L$ of $F$. We have $ L\cap
q^{-1}(\Mthick(\epsilon))\ne\emptyset$, i.e. 
$ q(L)\cap
\Mthick(\epsilon)\ne\emptyset$. 
Since $S$ is $\epsilon$-\good, and $q(L)$ is a union of cells of $K_S$,
Condition (2) of Definition \ref{good-def} guarantees that $ q(L)\cap
\inter\Mthick(\epsilon)\ne\emptyset$, i.e. 
$ L\cap
q^{-1}(\inter\Mthick(\epsilon))\ne\emptyset$. 
Note also that since $L$ is a face of
$F$, we have $L\in\realcalX_\tS$.

According to Condition (1) of  Definition \ref{good-def},
$\realcalX_\tS$ is weakly simple, and hence $L$ has valence $3$. As
$X_0$ and $X_1$ are distinct three-dimensional elements of $\realcalX_\tS$
having $L$ as a face,
there is exactly one three-dimensional element $X$ of $\realcalX_\tS$
which is distinct from $X_0$ and $X_1$ and has $L$ as a
face. Furthermore,  according to Conclusion (3) of
Proposition \ref{about 3}, $F_i:=X\cap X_i$
is a common two-dimensional face of $X$ and  $X_i$ for $i=0,1$.

According to the definitions given above, the intersections $F\cap
F_0$, $F\cap F_1$ and $F_0\cap F_1$ are all equal to $X\cap X_0\cap
X_1$. But Conclusion (2) of Proposition \ref{about 3} gives
$X\cap X_0\cap
X_1=L$, and hence
$F\cap
F_0=F\cap F_1=F_0\cap F_1=L$.

For $i=0,1$, since  the open set $q^{-1}(\inter\Mthick(\epsilon))$
meets $L\subset F_i$, it must meet the interior of $F_i$. Hence $\Phi_i^{-1}(\inter\Mthick(\epsilon))$
meets the interior of the face $V_i:=J_i(F_i)$ of $\calD_i$.
Since $\redcalT $ is $\epsilon$-\adapted, it follows that, for $i=0,1$, the
interior of $V_i $ contains a
(unique) point $\wouldbeV_i$ of $\Phi_i^{-1}(\redcalT )$.
We now have $(H_i,\wouldbeV_i)\in\scrP_{S,\redcalT }$ for $i=0,1$. 
We set
$\theta_i=\eta_{H_i,\wouldbeV_i}$ for $i=0,1.$ We set
$\eta_0=\theta_{0}$ and
$\eta_1=\bartheta_{1}$.

For $i=0,1$ we have
$\wasC_i\cap V_i
=J_i(F)\cap J_i(F_i)=J(L)$,
 so that 
$\wasC_i\cap V_i$ is a one-dimensional face of
$\calD_i$.
But we have $\wasC_i^\eta=\wasC_i$ and $\wasC_i^{\eta_i}=\wasC_0^{\theta_i}=V_i$. Thus $\wasC_i^{\eta}\cap 
\wasC_i^{\eta_i}$ is a one-dimensional face of
$\calD_i$. For $i=0$, this says that the oriented edges $\eta$ and $\eta_0$ of $G^{S,\redcalT }$, which
both have initial vertex $p_0$, are initially adjacent; and for $i=1$, it says that the oriented edges $\eta$ and $\eta_1$, which
both have terminal vertex $p_1$, are terminally adjacent.  This
establishes Assertions (1) and (2) of the conclusion of the lemma.

Now suppose that for each $i\in\{0,1\}$ we have chosen a
parametrization $\alpha_i$ of the oriented edge $\eta_i$. Set
$\beta_0=\alpha_0$ and $\beta_1=\baralpha_1$, so that $\beta_i$ is a
parametrization of 
the oriented edge
$\theta_i$ for $i=0,1$. 
We have
$\beta_i(0)=p_{i}$, and hence 
there is a lift $\tbeta_i$ of $\beta_i$ to $\HH^3$ such that
$\tbeta_i(0)=P_i$. 
Set $R_i=\tbeta_i(1)$ and $Y_i=X_{R_i}\in\realcalX_\tS$.

Let an index $j\in\{0,1\}$ be
given. Since $\beta_j$ is a
parametrization of 
the oriented edge
$\theta_j$, and since the lift $\tbeta_j$ of $\beta_j$ has initial
point $P_j$ and terminal point $R_j$, we may apply 
Lemma
\ref{match 'em}, with 
$\theta_j$
playing the role of $\eta$, and with $P_j$ and $R_j$ playing the
respective roles of $P_0$ and $P_1$. Then $X_j$ plays the role of
$X_{0}$ in 
Lemma
\ref{match 'em}, while $Y_j$ plays the role of $X_{1}$. Furthermore, the roles of
$\Hinit^\eta$, $\calD_{\Hinit^\eta,S}$, $J_{X_0,S}$,  and $U_0^\eta$ are played by
$H_j$, $\calD_j$, $J_j$,  and $V_j$ respectively. It now follows from
Lemma
\ref{match 'em} that
$q(\inter X_j)=H_j^\eta$, that $F_j':=X_j\cap Y_j$ is a common
two-dimensional face of $X_j$ and $Y_j$, and that
$J_j(F_j')=V_j$. Since $J_j:X_j\to\calD_j$ is a bijection and
$J_j(F_j')=V_j$, it follows that $F_j'=F_j$, i.e. that $X_j\cap
Y_j=F_j$. But we have $X_j\cap X=F_j$, and since $\realcalX_\tS$ is a
polyhedral complex whose underlying space $\HH^3$ is a $3$-manifold,
there can be only one three-dimensional element of $\realcalX_\tS$ whose
intersection with $X_j$ is the two-dimensional face $F_j$ of $X_j$.
Hence $Y_j=X$.

As the equality $Y_j=X$  holds for $j=0$ and for $j=1$, we have
$Y_0=Y_1$. Since $R_j$ is the unique point of $\tS$ in $Y_j$ for each $j\in\{0,1\}$, we have
$R_0=R_1$, i.e. $\tbeta_0(1)=\tbeta_1(1)$.

In particular we have $\beta_0(1)=\beta_1(1)$,
i.e. $\alpha_0(1)=\alpha_1(0)$; this gives Assertion (3) of the
lemma. Furthermore, since $\tbeta_0(1)=\tbeta_1(1)$, there is a
well-defined path $\talpha':=\tbeta_0\star\overline{\tbeta_1}$ in
  $\HH^3$. The path $\talpha'$ is a lift of $\alpha_0\star\alpha_1$,
  and, like $\talpha$, it has initial point $P_0$ and terminal point
  $P_1$. Since the covering space $\HH^3$ of $M$ is simply connected,
  it follows that $\alpha_0\star\alpha_1$ and $\alpha$ are
  homotopic. This is Assertion (4).
\EndProof

We conclude this section with the following additional fact about a
graph of the form $G^{S,\redcalT }$, which will be needed in
the proof of Proposition \ref{betwixt and between}.

\Lemma\label{they're short}
Let $M$ be a hyperbolic $3$-manifold of finite volume, let $\epsilon$
be a strict Margulis number for $S$,
let $S$ be a maximal $\epsilon$-thick $\epsilon$-\net\ in $M$, and let
$\redcalT $ be an $\epsilon$-\adapted\ \dotsystem\ for $S$. Then each edge of
the graph
$\redG^{S,\redcalT }$ is a rectifiable open arc of length less than $2\epsilon$.
\EndLemma


\Proof
Let $e$ be an edge of $G^{S,\redcalT }$. According to the material reviewed in
\ref {graph review}, 
there exist a point $\tau\in \redcalT $ and elements $(H\redsubnought ,\wasU\redsubnought )$ and $(H\redsubone ,\wasU\redsubone )$
of $
\redscrPST
$,  
interchanged by the involution
 $\sigma_{S,\redcalT }$, such that $\Phi_{H\redsubnought
   ,S}(\wasU\redsubnought )=\Phi_{H\redsubone ,S}(\wasU\redsubone
 )=\tau$, 
and $e$ is the
union of $\{\tau\}$ with the half-open  arcs
$B\redsubnought :=\{\tau\}
\cup\Phi_{H\redsubnought ,S}(\inter\ell_{H\redsubnought ,\wasU\redsubnought })$ and $B\redsubone :=\{\tau\}
\cup\Phi_{H\redsubone ,S}(\inter\ell_{H\redsubone ,\wasU\redsubone })$, whose intersection is their
common endpoint $\tau$. 

Now let us write $M=\HH^3/\Gamma$, where $\Gamma\le\isomplus(\HH^3)$
is discrete and torsion-free and has cofinite volume, let $q:\HH^3\to
M$ denote the quotient map, set $\tS=q^{-1}(S)$, and for 
$i=0,1$, 
let us fix a component
$\tH_i$ of $q^{-1}(H_i)$. Then for 
$i=0,1$ 
we have $\tH_i=\inter
X_{P_i}\in\realcalX_\tS$ 
for some $P_i\in\tS$. 
The closure of $B_i$ is
a closed arc which admits a lift to $X_{P_i}$; this lift is a geodesic
arc joining $P_i$ to some point $\ttau_i\in q^{-1}(\tau)$. Hence $B_i$
is a half-open geodesic arc, and $\length B_i=\dist(P_i,\ttau_i)$.

Since $\redcalT $ is $\epsilon$-adapted, Definition \ref{one way} gives
$\tau\in \redcalT \subset\Mthick(\epsilon)$. For 
$i=0,1$ 
we therefore have $\ttau_i\in X_{P_i}\cap q^{-1}(\{\tau\})\subset
X_{P_i}\cap q^{-1}(\Mthick(\epsilon))$. 
But since $S$ is a maximal $\maybeepsilon$-thick $\maybeepsilon$-\net\ for
$M$, it follows from Conclusion (5)
of \cite[Lemma 4.8]{ratioI} that
for every $P\in\tS$ we have $X_P\cap
  q^{-1}(\Mthick(\maybeepsilon))\subset\nbhd_{\maybeepsilon}(P)$.
Hence
  for 
$i=0,1$ 
we have $\ttau_i\in\nbhd_\epsilon(P_i)$, so that 
$\length B_i=\dist(P_i,\ttau_i)<\epsilon$. It now follows that
$\length e=\length B\redsubnought  +\length B\redsubone <2\epsilon$.
\EndProof

\section{Coloring systems}
\DefinitionRemark\label{coloring system def}
Let $S$ be a non-empty finite subset of an orientable hyperbolic
$3$-manifold $M$, and let $\redcalT $ be a \dotsystem\ for $S$. We recall from
\ref{graph review} that $G^{S,\redcalT }$ is a graph whose underlying space is
$\scrG^{S,\redcalT }\subset M$, and that
$\veccale^{S,\redcalT }$ denotes the set of all oriented edges of the graph
$\redG^{S,\redcalT }$.
We define a {\it coloring system} for $(S,\redcalT )$ in $M$ to
be a map of sets $\frakc :\veccale^{S,\redcalT }\to\{1,2,3,4\}$ 
satisfying the following condition:
\Claim\label{not the same}
If $\eta$ and $\eta'$ are elements of
$\veccale^{S,\redcalT }$ such that $\init_\eta=\init_{\eta'}$, and $\eta$ and $\eta'$
are initially adjacent
(see \ref{adjacency def}), 
then $\frakc (\eta)\ne \frakc(\eta')$.
\EndClaim
Note that if $\frakc $ is a coloring system for $(S,\redcalT )$ in $M$, and if $\eta$ and $\eta'$ are elements of
$\veccale^{S,\redcalT }$ such that $\term_\eta=\term_{\eta'}$, and $\eta$ and $\eta'$
are terminally adjacent 
(see \ref{adjacency def}), 
then $\frakc (\overline{\eta})\ne \frakc(\overline{\eta'})$. This follows upon applying \ref{not the same} with $\overline{\eta}$ and $\overline{\eta'}$ playing the respective roles of 
$\eta$ and $\eta'$.
\EndDefinitionRemark

\Proposition\label{because four colors}
If $S$ is any non-empty finite subset of an orientable hyperbolic
$3$-manifold $M$, and  $\redcalT $ is any \dotsystem\ for $S$, then there exists a coloring system for $(S,\redcalT )$ in $M$.
\EndProposition

\Proof
For each $3$-cell $H$ of $K_S$, let $\calf_H$
denote the set of all two-dimensional faces of $\maybecalD_{H,S}$, and
choose a function $\kappa_H $ with domain $\calf_H$, taking values in the
set $\{1,2,3,4\}$, and having the following
property:
\Claim\label{four suffice}
For any two two-dimensional faces $F$ and $F'$  of $\maybecalD_{H,S}$
such that $F\cap F'$ is a one-dimensional polyhedron, we have
$\kappa_H (F)\ne\kappa_H (F')$.
\EndClaim

In the special case where $\calf_H$ is finite, the existence of a
function $\kappa_H $ for which \ref{four suffice} holds is the most
classical statement of the Four-Color Theorem, which was proved in
\cite{haken-appel} and \cite{rsst}.  
In the general case, where $\calf_H$ may be infinite, the existence of
such a function follows from the Four-Color Theorem for possibly
infinite planar graphs, which is deduced 
from  the finite version by the de Bruijn-Erd\"os theorem
\cite{infinite-coloring}. The latter theorem states that for any
positive integer $k$, a graph has a
$k$-coloring if each of its finite subgraphs has a $k$-coloring.

Now define a map  of sets $\frakc :\veccale^{S,\redcalT }\to\{1,2,3,4\}$ by
setting
$$
\frakc(\eta)=\kappa_{\Hinit^\eta}(\Cinit^\eta).
$$
(The notation $\Hinit^\eta$ and
$\Cinit^\eta$ is defined in \ref{before adjacency}.)

In order to show  that $\frakc$ is a coloring system, we must verify
Condition \ref{not the same}. Suppose that  $\eta$ and
$\eta'$ are oriented edges of $\redG^{S,\redcalT }$ such that
$\init_\eta=\init_{\eta'}$, and that $\eta$ and $\eta'$ are initially
adjacent. Set $H:=\Hinit^\eta$.
Since $\init_\eta=\init_{\eta'}$, we  have $H=\Hinit^{\eta'}$, and 
   $\Cinit^\eta$ and $\Cinit^{\eta'}$ are two-dimensional faces of
   $\maybecalD_{H,S}$. By the definition of initial adjacency given in
   \ref{adjacency def},
   $\Cinit^\eta\cap\Cinit^{\eta'}$ is a one-dimensional face of $\maybecalD_{H,S}$.
It then follows from \ref{four suffice} that
$\kappa_H(\Cinit^\eta)\ne\kappa_H(\Cinit^{\eta'})$, i.e. $\frakc(\eta)\ne\frakc(\eta')$.
\EndProof

\Definition\label{coloring assignment def}
Let $S$ be a finite subset of an orientable hyperbolic $3$-manifold $M$. We define a {\it color assignment} for $S$ in $M$ to be a map of sets $\fraka :S\to\{1,2,3,4\}$. 
\EndDefinition

\Definition\label{distinguished def}
Let $S$ be a finite subset of an orientable hyperbolic $3$-manifold $M$, let $\fraka $ be a color assignment for $S$, let $\redcalT $ be a \dotsystem\ for $S$, and let $\frakc $ be a
coloring system for $(S,\redcalT )$. We will say that an oriented edge $\eta$
of $\redG^{S,\redcalT }$ is {\it $(\frakc,\fraka)$-distinguished} if $\frakc
(\eta)=\fraka(\init_\eta)$. We will say that an
edge $e$ of $\redG^{S,\redcalT }$ is {\it $(\frakc,\fraka)$-doubly distinguished} if both the orientations of $e$ are $(\frakc,\fraka)$-distinguished oriented edges of $\redG^{S,\redcalT }$.
\EndDefinition

\Proposition\label{one-sixteenth}
Let $S$ be a finite subset of an orientable hyperbolic $3$-manifold $M$, let $\redcalT $ be a \dotsystem\ for $S$,  let $\frakc $ be a
coloring system for $(S,\redcalT )$, and let $\cale_0$ be a set of  (unoriented) edges of $\redG^{S,\redcalT }$ such that no loop belongs to $\cale_0$. For each color assignment $\fraka $ for $S$,  let $\cale_0^\fraka\subset\cale_0$ denote the set of all $(\frakc,\fraka)$-doubly distinguished edges in $\cale_0$. Then there exists a color assignment $\fraka_0 $ for $S$ such that
$$\#\cale_0^{\fraka_0}\ge\frac{\#\cale_0}{16}.$$
\EndProposition

\Proof
The vertex set of
$\redG^{S,\redcalT }$ is $S$ (see \ref{graph review}).
We set $s=\#(S)$.

Let $\cala$ denote the set of all color
assignments for $S$. Then $\#(\cala)=4^s$.

For each $e\in\cale_0$, let $\cala_e$ denote the set of all color
assignments $\fraka\in\cala$  such that $e$ is $(\frakc,\fraka)$-doubly distinguished. Then we
have
\Equation\label{switch}
\sum_{\fraka\in\cala}\#(\cale_0^\fraka)=\sum_{e\in\cale_0}\#(\cala_e).
\EndEquation

Suppose that $e$ is an arbitrary element of $\cale_0$. Let
$\eta^{(1)},\eta^{(2)}\in \veccale^{S,\redcalT }$ denote the two orientations
of $e$ (so that $\eta^{(2)}=\overline{\eta^{(1)}}$),
and set $v_i=\init_{\eta^{(i)}}$ for $i=1,2$. (Thus we have $v_1=\term_{\eta^{(2)}}$ and
$v_2=\term_{\eta^{(1)}}$.) By hypothesis
$e$ is not a loop, and hence $v_1\ne v_2$.
For $i=1,2$, according to the definitions, $\eta^{(i)}$ is $(\frakc,\fraka)$-distinguished
for a given 
color
assignment $\fraka $ if and only if $\frakc (\eta^{(i)})=\fraka(v_i)$. Thus we
have $\fraka\in \cala_e$ if and only if
$\fraka (v_1) =\frakc(\eta^{(1)})$ and $\fraka (v_2)=\frakc(\eta^{(2)})$.

Let $\calF_{e}$ denote the set of all maps of sets
from $S-\{v_1,v_2\}$ to $\{1,2,3,4\}$; thus $\#(\calF_e)=4^{s-2}$. If
$f\in F_{e}$ is given, then since $v_1\ne v_2$, there is a
well-defined
color
assignment $\fraka $ for $S$ given by setting $\fraka (v)=f(v)$ for
$v\in S-\{v_1,v_2\}$, and $\fraka (v_i)=\frakc(\eta^{(i)})$ for $i=1,2$. This is
the unique assignment $\fraka \in\cala_e$ which extends the given map $f\in\calF_e$. This
shows that the restriction map $\cala_e\to\calF_e$ is a bijection, and
hence that $\#(\cala_e)=\#(\calF_e)=4^{s-2}$. As this equality holds
for every $e\in\cale_0$, (\ref{switch}) becomes 
\Equation\label{new switch}
\sum_{\fraka\in\cala}\#(\cale_0^\fraka)=\sum_{e\in\cale_0}4^{s-2}=4^{s-2}\cdot\#(\cale_0).
\EndEquation
Since $\#(\cala)=4^s$, at least one term on the left-hand side of
(\ref{new switch}) must be greater than or equal to
$4^{s-2}\cdot\#(\cale_0)/4^s=\#(\cale_0)/16$, and the assertion follows.
\EndProof

\Proposition\label{homotoping prop}
Let $M$ be  a
finite-volume 
orientable hyperbolic $3$-manifold, 
let $\epsilon$ be a strict Margulis number for $M$, and let $S\subset \Mthick(\epsilon)$ be a non-empty finite set which is $\epsilon$-\good.
Let $\redcalT $ be a \dotsystem\ for $S$
which is $\epsilon$-\adapted.
Let $\mathfrak a$ be a color assignment for $S$, 
and let $\mathfrak c$ be a
coloring system for $(S,\redcalT )$.
Let $\redGndd$ denote the subgraph of $\redG^{S,\redcalT }$ consisting of all
vertices of $\redG^{S,\redcalT }$ and all edges of  $\redG^{S,\redcalT }$ that are
not $(\frakc,\fraka)$-doubly distinguished. 
Then every path in 
$\scrG^{S,\redcalT }=|G^{S,\redcalT }|$
whose endpoints are vertices is fixed-endpoint homotopic in $M$ to some path in $|\redGndd|$.
\EndProposition

\Proof
Any path in $\scrG^{S,\redcalT }$
whose endpoints are vertices is a composition of finitely many paths
each of which is a parametrization of some  edge of
$\redG^{S,\redcalT }$. Hence it suffices to prove the proposition in the
special case of a path $\alpha$ which is itself a parametrization of a
single edge $e$ of $\redG^{S,\redcalT }$.

If $e$ is an edge of $\redGndd$, there is nothing to prove. For the
rest of the proof we will assume that $e$ is not an edge of $\redGndd$,
i.e. that $e$ is $(\frakc,\fraka)$-doubly distinguished. 

Let $\eta$ denote the orientation of $e$ determined by the
parametrization $\alpha$. Thus $\eta$ and $\bareta$ are both
$(\frakc,\fraka)$-distinguished. We set $v\redsubnought=\init_\eta$ and $v\redsubone=\term_\eta$.

Since $S$ is $\epsilon$-\good\ and $\redcalT $ is
$\epsilon$-adapted, Lemma \ref{and what if it is} gives oriented edges
$\eta\redsubnought$ and $\eta\redsubone$ of $\redG^{S,\redcalT }$ such that Conditions 
(1)--(4)
of that lemma hold. For 
$i=0,1$ 
we denote by $e_i$ the underlying edge
of $\eta_i$,

According to 
Condition (1)
of Lemma \ref{and what if it is}, we have
$\init_{\eta\redsubnought}=v\redsubnought$, and $\eta\redsubnought$ is initially adjacent to $\eta$. By
the definition of a coloring system it follows that $\frakc (\eta\redsubnought)\ne
\frakc(\eta)$. But we have $\frakc (\eta)=\fraka(v\redsubnought)$ since $\eta$ is $(\frakc,\fraka)$-distinguished,
and therefore $\frakc (\eta\redsubnought)\ne \fraka(v\redsubnought)$. This means that $\eta\redsubnought$ is not
$(\frakc,\fraka)$-distinguished; in particular, $e\redsubnought$ is not $(\frakc,\fraka)$-doubly
distinguished, i.e. $e\redsubnought$ is an edge of $\redGndd$.

According to 
Condition (2)
of Lemma \ref{and what if it is}, we have 
$\term_{\eta\redsubone}=v\redsubone$, and $\eta\redsubone$ is terminally adjacent to $\eta$. It
was observed in \ref{coloring system def} that this implies that $\frakc (\bareta\redsubone)\ne
\frakc(\bareta)$. But we have $\frakc (\bareta)=\fraka(v\redsubone)$ since $\bareta$ is $(\frakc,\fraka)$-distinguished,
and therefore $\frakc (\bareta\redsubone)\ne \fraka(v\redsubone)$. This means that $\bareta\redsubone$ is not
$(\frakc,\fraka)$-distinguished; in particular, $e\redsubone$ is not $(\frakc,\fraka)$-doubly
distinguished, i.e. $e\redsubone$ is an edge of $\redGndd$.

According to Conditions (3) and (4) of Lemma \ref{and what if it is},
$\alpha$ is homotopic to a path $\alpha'$ lying in the union of the closures of
the edges $e\redsubnought$ and $e\redsubone$. Since we have seen that $e\redsubnought$ and $e\redsubone$ are
edges of $\redGndd$, the path $\alpha'$ lies in $|\redGndd|$.
\EndProof

\section{
Conditional bounds on the rank of the fundamental group
}

 The bounds  given in this section for the rank of $\pi_1(M)$, where
 $M$ is a  hyperbolic $3$-manifold,
 are ``conditional'' in the sense that, in addition to depending on
 the volume of $M$ (or related quantities), they depend on
 bounds for the ranks of certain special subgroups of $\pi_1(M)$. The
 needed bounds for the ranks of these special subgroups will be
 supplied in the next section, and in the sequel to this paper.

\Proposition\label{new what's new}
Let $M$ be a finite-volume orientable hyperbolic $3$-manifold, and let
$\epsilon$ be a strict Margulis number for $M$. 
Let 
$S$ be a finite subset of $M$ which is
$\epsilon$-good,
and let 
$\redcalT $ be an $\epsilon$-adapted \dotsystem\ for $S$
(so that $S$ is the vertex set of the graph 
 $\redG^{S,\redcalT }$ by \ref{graph review}).
For each $v\in\redS$,
and each loop $\ell$ of $\redG^{S,\redcalT }$ based at the vertex $v$, let $[[\ell]]$
denote the element of $\pi_1(M,v)$, 
well-defined up to inversion,
which is the fixed-endpoint homotopy class of a parametrization of
$\ell$. For each  $v\in\redS$, suppose that $\wasbarX_v$
is a
 (possibly trivial) subgroup of $\pi_1(M,v)$ which contains
the element $[[\ell]]$ for every loop $\ell$ 
of $\redG^{S,\redcalT }$ based at $v$.

Then
\begin{enumerate}
\item
the graph 
$\redG^{S,\redcalT}$ is connected, and 
the inclusion homomorphism $\pi_1(\redG^{S,\redcalT})\to\pi_1(M)$ is
  surjective.
\end{enumerate}
 Furthermore,
there is a 
connected 
subgraph $\redG^\dagger$ of $\redG^{S,\redcalT }$ such that the
following  conclusions
 (2)--(4) hold.

\begin{enumerate}
\item [(2)]
The 
vertex set of $\redG^\dagger$ is all of $S$.
\item[(3)]
If $p$ is a point of $S$,  if for each $v\in\redS$  we are given a
path $\beta_v$ in $|\redG^\dagger|$ from $p$ to $v$, 
and if we let
$\wasbarY_v\le \pi_1(M,p)$ denote the preimage of $\wasbarX_v$ under the
isomorphism from $\pi_1(M,p)$ to $\pi_1(M,v)$ defined by $\beta_v$,
then $\pi_1(M,p)$ is generated by the family of subgroups
consisting of the image of the inclusion homomorphism
$\pi_1(|\redG^\dagger|,p)\to\pi_1(M,p)$ 
and the subgroups of the form $\wasbarY_v$, where
$v$ ranges over  
$\redS$.
\item[4)]
If $E$, $L$ and $\reds=\#(\redS)$ denote respectively the number of edges, loops and
vertices of $\redGST$, then the first betti number of $|\redG^\dagger|$
is at most 
$$
\frac{15}{16}(E-L)-\reds+1.
$$
\end{enumerate}
Finally:
\begin{enumerate}
\item[(5)]
If $E$ and $\reds$ are defined as in 
(4), 
we have 
$$
\rank\pi_1(M)\le\frac{15}{16}E-\reds+1+\frac1{16}\sum_{v\in\redS}\rank\wasbarX_v.
$$
\end{enumerate}
\EndProposition


\Proof
Set $\Theta=\Mthick(\epsilon)$. Since $\epsilon$ is a strict Margulis number for $M$, it follows from \cite[Proposition 4.5]{ratioI} that $\Theta$ is a
connected (and in particular non-empty)
$3$-manifold-with-boundary, 
and that the inclusion homomorphism $\pi_1(\Theta)\to\pi_1(M)$
is surjective. Since $\Theta\ne\emptyset$, the maximality of $S$ implies that $S\ne\emptyset$. 


It follows from the definition of  an $\epsilon$-\adapted\ \dotsystem\
for $S$ (see \ref{one way})
that
for every  $2$-cell $C$ 
of $K_S$
with $C\cap\Theta\ne\emptyset$, we have $\redcalT \cap
C\cap\Theta\ne\emptyset$. 
This shows that $M$, $S$, $\Theta$ and $\redcalT $ satisfy the hypotheses of \cite[Proposition
2.14]{ratioI}. 
Now Conclusion (1) of \cite[Proposition
2.14]{ratioI} asserts that 
any two points of $S$ that lie in the same component of
$\Theta$ also lie in the same component of the set
$\scrG^{S,\redcalT }=|G^{S,\redcalT }|$.
Since $\Theta$ is connected, it follows that $\scrG^{S,\redcalT }$ is
connected. Furthermore, given the connectedness of $\scrG^{S,\redcalT
}$ and $\Theta$,
Conclusion (2) of \cite[Proposition
2.14]{ratioI} asserts that for any base point $p\in S$ we have 
$\image(\pi_1(\Theta,p)\to\pi_1(M,\willbep))\le \image(\pi_1(\scrGST,\willbep)\to\pi_1(M,\willbep)),$
where the unlabeled arrows denote inclusion
homomorphisms.
Since 
the inclusion homomorphism $\pi_1(\Theta)\to\pi_1(M)$
is surjective, we deduce that
the inclusion homomorphism $\pi_1(\scrGST)\to\pi_1(M)$ is
surjective. This establishes Conclusion (1) of the present proposition.

According to Proposition \ref{because four colors}, we may fix a
coloring system $\frakc $ for $(S,\redcalT )$. We denote by $\cale_{\rm nl}$ the set
of all edges of $\redGST$  that are not loops, and we apply
Proposition \ref{one-sixteenth} with $\cale_{\rm nl}$ playing the role
of $\cale_0$ in that proposition. This gives a color assignment $\fraka $
for $S$ (denoted $\fraka_0$ in the conclusion of Proposition
\ref{one-sixteenth}) such that at least $1/16$ of the edges in $\cale_{\rm nl}$ are
$(\frakc,\fraka)$-doubly distinguished.  Thus if we denote by $\redcalestar$ the set
of all edges of $\redGST$ that are not $(\frakc,\fraka)$-doubly distinguished,
we have
\Equation\label{java joe}
\#(\calenl\cap\redcalestar)\le\frac{15}{16}\cdot\#(\calenl).
\EndEquation

As in the statement of the present proposition, we 
will denote by
$E$ and $\reds =\#(\redS)$ the numbers of edges and vertices,
respectively, of $\redGST$; and 
we will denote by $L$ the number of loops in
$\redGST$.
We then have $\#(\calenl)=E-L$. Thus (\ref{java joe}) may be
rewritten in the form
\Equation\label{more java}
\#(\calenl\cap\redcalestar)\le\frac{15}{16}(E-L).
\EndEquation

As in 
the statement of Proposition \ref{homotoping prop}, we denote by $\redGndd$  the subgraph of $\redG^{S,\redcalT }$ consisting of all
vertices of $\redG^{S,\redcalT }$, and all edges of  $\redG^{S,\redcalT }$ that are
not $(\frakc,\fraka)$-doubly distinguished; thus $\redGndd$ has vertex
set $S$ and edge set $\redcalestar$.
According to 
Proposition \ref{homotoping prop}, every path in $|\redGST|$
whose endpoints are vertices is fixed-endpoint homotopic in $M$ to
some path in $|\redGstar|$.  
Combining this with Conclusion (1) of the
present proposition, which has been established, 
we
deduce:

\Claim\label{dinkum beano}
The graph $\redGstar$ is connected, and the inclusion homomorphism
$\pi_1(|\redGstar|)\to\pi_1(M)$ is surjective.
\EndClaim

Let $\redG^\dagger$ 
denote the subgraph of $\redGstar$ consisting of all vertices 
in $\redS$
and all edges of $\redGstar$ that are not loops;
thus $\redG^\dagger$ has vertex set $S$ and edge set
$\calenl\cap\redcalestar$. 
Since
$\redGstar$ is connected, $\redG^\dagger$ is also connected. 
Note
also that, according to the definition of
$\redG^\dagger$, 
Conclusion (2) of the present proposition holds.

 We 
now fix, arbitrarily, a base point $p\in S$, and 
denote by 
$Z\le\pi_1(|\redGstar|,p)$ 
the image of the inclusion homomorphism $
\pi_1(|\redG^\dagger|,
p
)\to\pi_1(|\redGstar|,
p
)$.

For each 
 $v\in\redS$ 
we denote by $\scrL_v$ the connected subgraph of
 $\redGstar$ whose only vertex is $v$, and whose edges are all the
 loops in $\redGstar$ based at $v$. (Of course for a given 
$v\in\redS$ 
we may well have $|\scrL_v|=\{v\}$.) 
Suppose that for each  
$v\in\redS$ 
we are given a path $\beta_v$ in $|\redG^\dagger|$ from $p$ to $v$, 
and let  $i_v:\pi_1(|\redGstar|,p)\to\pi_1(|\redGstar|,v)$  denote
the isomorphism defined by the path $\beta_v$. 
For each $v$ we denote by $X_v\le\pi_1(|\redGstar|,v)$ the image of the inclusion homomorphism 
$\pi_1(|\scrL_v|,v)\to\pi_1(|\redGstar|,v)$, and we set
$Y_v=i_v^{-1}(X_v)\le\pi_1(|\redGstar|,p)$. Then:
\Claim\label{it's clear}
The group $\pi_1(|\redGstar|,p)$ is
generated by the family of subgroups consisting of $Z$ and the
subgroups of the form $Y_v$, where $v$ ranges over  
$\redS$.
\EndClaim

Now if
$j:\pi_1(|\redGstar|,p)\to\pi_1(M,p)$ denotes the inclusion
homomorphism, then 
 the 
image of the inclusion homomorphism
$\pi_1(|\redG^\dagger|),p)\to\pi_1(M,p)$ is equal to
$j(Z)$. 
On the other hand, the hypothesis that 
$\wasbarX_v$
contains $[[\ell]]$ for every loop $\ell$
of $\redG^{S,\redcalT }$ based at $v$,
implies that
for each $v\in\redS$,
the group denoted by  $\wasbarY_v$ in the statement of the
proposition contains $j(Y_v)$. 
Furthermore, according to \ref{dinkum beano} we have 
$j(\pi_1(|\redGstar|,p))=\pi_1(M,p)$. 
Hence Conclusion (3) of
the present proposition follows from \ref{it's clear} upon
applying the homomorphism $j$.



To show that  Conclusion (4) holds, 
let $\wasbeta$ denote the
first betti number of the graph
$\redG^\dagger$. 
The vertex set of $\redG^\dagger$ is by definition
equal to 
$\redS$, which has cardinality $\reds $. 
Since the edge set of
$\redG^\dagger$ is $\calenl\cap\redcalestar$,
the Euler characteristic of
$|\redG^\dagger|$ is equal to $\reds -\#(\calenl\cap\redcalestar)$;
thus we have
\Equation\label{betti boop}
\wasbeta=\#(\calenl\cap\redcalestar)-\reds +1.
\EndEquation
Combining (\ref{betti boop}) with 
(\ref{more java}), 
we obtain
$$
\wasbeta\le\frac{15}{16}(E-L)-\reds +1,
$$
which is Conclusion (4).

It remains to show that 
Conclusion (5) holds. 
For this
purpose, we first note that according to Conclusion (3), $\pi_1(M)$ is
generated by a family of subgroups, of which one is a homomorphic
image of the group $\pi_1(|\redG^\dagger|,p)$, while the rest are the
subgroups of the form $\wasbarY_v$ for $v\in\redS$. By definition
$\wasbarY_v$ is isomorphic to $\wasbarX_v$. Hence we have 
\Equation\label{some numbers}
\rank\pi_1(M)\le\rank \pi_1(|\redG^\dagger|)+
\sum_{v\in\redS}\rank\wasbarX_v.
\EndEquation
Since $\rank \pi_1(|\redG^\dagger|)$ is equal to the first betti number
of $|\redG^\dagger|$, we can combine (\ref{some numbers}) with Conclusion
(4) to deduce that 
\Equation\label{where was it}
\rank\pi_1(M)\le\frac{15}{16}(E-L)-\reds +1+
\sum_{v\in\redS}\rank\wasbarX_v.
\EndEquation

On the other hand, for each  
$v\in\redS$, 
the
group $\wasbarY_v$ is a homomorphic image of $Y_v$; thus 
$\rank A_v=
\rank \wasbarY_v$ 
is
bounded above by $\rank Y_v$. The group $Y_v$ is isomorphic to $X_v$,
so that its rank is equal to the first betti number of $\scrL_v$,
which is the number of loops in $\redGstar$ based at $v$. As $\redGstar$ 
is a subgraph of $\redGST$, the rank of $Y_v$ is in particular bounded
above by the number of loops in $\redGST$ based at $v$. It follows
that
\Equation\label{rusty racoon}
\sum_{v\in\redS}\rank A_v\le L.
\EndEquation

By (\ref{where was it}) and
(\ref{rusty racoon}), we have
$$\begin{aligned}
\rank\pi_1(M) &\le\frac{15}{16}(E-L)-\reds +1+
\frac{15}{16}\sum_{v\in\redS}\rank\wasbarX_v+\frac{1}{16}\sum_{v\in\redS}\rank\wasbarX_v\\
&\le\frac{15}{16}(E-L)-\reds +1+
\frac{15}{16}L+
\frac{1}{16}\sum_{v\in\redS}\rank\wasbarX_v
\\
&=\frac{15}{16}E-\reds+1+\frac1{16}\sum_{v\in\redS}\rank\wasbarX_v,
\end{aligned}
$$
which gives Conclusion (5).
\EndProof

The following result, Proposition \ref{betwixt and between}, involves
a number of functions and a set that are defined in \cite{ratioI}. As
in \cite[Subsection 3.1]{ratioI},  for
$\mayber>0$ we denote by
$B(\mayber)$ the function $\pi(\sinh(2\mayber)-2\mayber)$, which gives the volume of a ball  in $\HH^3$ of
radius $\mayber$; and as in \cite[Subsection 4.1]{ratioI}, we denote
by $b(r)$ the function, also defined for $r>0$, and given by
$b(r)=B(r)/\density(r)$, where $\density(r)$ is B\"or\"oczky's density
bound for sphere-packings in hyperbolic space, the definition of which is reviewed in
detail in \cite[Subsection 4.1]{ratioI}. The rather involved but
elementary definitions of the set $\scrV_0\subset\RR^3$, and of the
function $\phi$, whose domain is $\scrV_0$, are given in Subsection 3.9
of \cite{ratioI},  and depend on definitions given in Subsections 3.1,
3.2,
3.4 and 3.7 of that paper.

According to \cite[Lemma 4.9]{ratioI}, if $\tryepsilon$ and $R$ are
positive numbers such that $2\tryepsilon<R<5\tryepsilon/2$, then for
every number $D$ in the interval 
$[R/2-\tryepsilon/4,\tryepsilon]$, 
we have
$(R-D,\tryepsilon/2,D)\in\scrV_0$. This fact is needed as background
for the statement of 
Proposition \ref{betwixt and between} below.


\Proposition\label{betwixt and between}
Let $M$ be a \willbefinitevolume\  orientable hyperbolic $3$-manifold. 
Suppose that $\epsilon$ is a Margulis number for $M$, and
that $R$ is a real number with 
$2\epsilon<R<5\epsilon/2$
(so that  
for
every  
$D\in[R/2-\epsilon/4,\epsilon]$
we have 
$(R-D,\epsilon/2,D)\in\scrV_0$ by \cite[Lemma 4.9]{ratioI}, and hence
$\phi(R-D,
\epsilon/2,
D)$ is defined). 
Let $c$ be a positive number such that (a) 
$\phi(R-D,
\epsilon/2,
D)>c$ for every $D\in[R/2-\epsilon/4,\epsilon]$, 
and (b) $B(\epsilon/2)>c$.
Suppose that $\rho$ is a non-negative integer such that, for every
point $v\in M$, the subgroup of $\pi_1(M,v)$ generated by all elements
represented by closed paths of length less than $2\epsilon$ based at
$v$ has rank at most $\rho$. Then we have
$$\rank\pi_1(M)\le1+\bigg\lfloor\frac{\vol M}{b(\epsilon/2)}\bigg\rfloor\cdot
\max\bigg(0,
\frac{15}{32}\bigg\lfloor\frac{B(R)-b(\epsilon/2)}c\bigg\rfloor-1+\frac{\rho}{16}\bigg).$$
\EndProposition

\Proof
We will first give the proof under the following additional
assumptions:
\Claim\label{additional}The number $\epsilon$ is a strict Margulis
number for $M$, and $M$ contains a maximal $\epsilon$-thick
$\epsilon$-\net\ which is $\epsilon$-good.
\EndClaim

Suppose that the assumptions stated in \ref{additional} hold, and fix
an  $\epsilon$-thick
$\epsilon$-\net\ $S$ in $M$ which is $\epsilon$-good. According to
\ref{one way}, we may fix an $\epsilon$-\adapted\ \dotsystem\ $\redcalT $ for
$S$. For each 
$v\in S$, 
let $A_v$ denote
the subgroup of $\pi_1(M,v)$ generated by all elements
represented by closed paths of length less than $2\epsilon$ based at
$v$. 

According to Lemma \ref{they're short}, for every loop $\ell$ 
of $\redG^{S,\redcalT }$ based at $v$, a parametrization of $\ell$ has
length less than
$2\epsilon$; 
hence, for every vertex $v$ of $\redG^{S,\redcalT }$ and every loop $\ell$ 
of $\redG^{S,\redcalT }$ based at $v$, the  subgroup  $\wasbarX_v$
of $\pi_1(M,v)$ contains 
the element $[[\ell]]$ represented by a parametrization
of $\ell$ (and well defined up to inversion). 
Thus all the hypotheses of
Proposition \ref{new what's new} hold with the choices of $S$, $\redcalT $ and
$A_v$ made above.
Hence the
inequality in 
Conclusion (5) 
of  Proposition \ref{new
  what's new} 
holds when 
$E$, $L$ and $s=\#(S)$ denote respectively the number of edges, loops and
vertices of $\redGST$.
According to the hypothesis of the present proposition we have
$\rank A_v\le\rho$ for each $v\in S$. 
Hence the inequality in
Conclusion (5) 
of Proposition \ref{new what's new}
implies
\Equation\label{whaddya think it's four}
\rank\pi_1(M)\le\frac{15}{16}E-s+1+\frac\rho{16} s.
\EndEquation

Let $\calc$ denote the set of all $2$-cells of $K_S$ which meet
$\Mthick(\epsilon)$. Then according to 
\cite[Subsection 2.13]{ratioI} (cf. Subsection \ref{gpc review} above), 
for each  
$v\in S$ 
there is a $3$-cell $H$ of $K_S$ such that
the
valence of 
the vertex $v$ in the graph $\redGST$  
is
the 
number of two-dimensional
faces of $ \maybecalD_{H,S}$ whose interiors are mapped by
$\Phi_{H,S}$ onto $2$-cells that meet $\redcalT $.
Since $S$ is a maximal $\epsilon$-thick $\epsilon$-\net, it follows
from Conclusion (2) of \cite[Proposition
4.10]{ratioI} that for any $3$-cell $H$ of $K_S$, this number is
bounded above by
\Equation\label{upper}
\bigg \lfloor\frac{B(R)-\maybeb(\tryepsilon/2)}c
\bigg\rfloor.
\EndEquation

Thus (\ref{upper}) is an upper bound for the valence in
$\redGST$ of any vertex $v\in S$. 
Hence
$$E\le \frac12\bigg \lfloor\frac{B(R)-\maybeb(\tryepsilon/2)}c
\bigg\rfloor\cdot s,$$
which with (\ref{whaddya think it's four}) gives
\Equation\label{i'll bet}
\rank\pi_1(M)\le1+s\cdot\bigg(\frac{15}{32}\bigg\lfloor\frac{B(R)-b(\epsilon/2)}c\bigg\rfloor-1+\frac{\rho}{16}\bigg).
\EndEquation

If we now set 
$$m=\max \bigg( 0,\frac{15}{32}\bigg\lfloor\frac{B(R)-b(\epsilon/2)}c\bigg\rfloor-1+\frac{\rho}{16}\bigg),
$$
then (\ref{i'll bet}) implies in particular that
$\rank\pi_1(M)\le1+sm$. 
On the other hand, according to 
Conclusion (1) of \cite[Proposition
4.10]{ratioI}, 
the quantity
$s$ is 
bounded above by
$\lfloor\vol(M)/b(\epsilon/2)\rfloor$. Since $m\ge0$, it follows that 
$$\rank\pi_1(M)\le1+sm\le1+\bigg\lfloor\frac{\vol M}{b(\epsilon/2)}\bigg\rfloor\cdot
m,$$
which gives the conclusion of the proposition in this case.

We now turn to the general case. It follows from Lemma \ref{when is
  it}, applied with $\epsilon$ playing the role of $\epsilon_0$, that
$\epsilon$ is the limit of a sequence $(\epsilon_n)_{n\ge1}$, with
$0<\epsilon_n<\epsilon$ for each $n$ (so that each $\epsilon_n$ is a strict
Margulis number for $M$), and such that for each $n$ there exists
a maximal
$\epsilon_n$-thick $\epsilon_n$-\net\ which is $\epsilon_n$-\good. Since
$2\epsilon<R<5\epsilon/2$, we have $2\epsilon_n<R<5\epsilon_n/2$ for
all sufficiently large $n$. 
The property (b) of the number $c$ which is stated in the hypothesis
of the present proposition
says that
$B(\epsilon/2)>c$; 
hence for sufficiently large $n$ we have
$B(\epsilon_n/2)>c$, i.e. Property (b) still holds when $\epsilon$ is
replaced by $\epsilon_n$. 

We claim that for
all sufficiently large $n$,
the property (a) of the number $c$ which is stated in the hypothesis
of the present proposition
also continues to hold 
when $\epsilon$ is
replaced by $\epsilon_n$; that is, $\phi(R-D,
\epsilon_n/2,
D)>c$ for every $D\in[R/2-\epsilon_n/4,\epsilon_n]$. Assume that this
is false. Then after passing to a subsequence, we may assume that for
each $n$ there is a $D_n\in[R/2-\epsilon_n/4,\epsilon_n]$ such that $\phi(R-D_n,
\epsilon_n/2,
D_n)\le c$ (where $\phi(R-D_n,
\epsilon_n/2,
D_n)\le c$ is defined because $(R-D_n,
\epsilon_n/2,
D_n)\in\scrV_0$ by 
\cite[Lemma 4.9]{ratioI}).
In particular the sequence $(D_n)_{n\ge1}$ is bounded, and
after again passing to a subsequence we may assume it  converges to
some limit $D_\infty$. Then  $D_\infty\in[R/2-\epsilon/4,\epsilon]$, so
that $(R-D_\infty,
\epsilon/2,
D_\infty)\in\scrV_0$; and Property (a) of the number $c$, as stated in
the hypothesis, implies that 
 $\phi(R-D_\infty,
\epsilon/2,
D_\infty)>c$. But since 
$\phi(R-D_n,
\epsilon_n/2,
D_n)\le c$ for each $n$, the continuity of $\phi$ on its domain
$\scrV_0$ (see
\cite[Subsection 3.9]{ratioI}) implies that 
 $\phi(R-D_\infty,
\epsilon/2,
D_\infty)\le c$. This contradiction establishes our claim.

Thus, for sufficiently large $n$, all the hypotheses of the present
proposition continue to hold when $\epsilon$ is replaced by
$\epsilon_n$. Since 
each
$\epsilon_n$ is a strict
Margulis number for $M$, and since for each $n$ there exists
a maximal
$\epsilon_n$-thick $\epsilon_n$-\net\ which is $\epsilon_n$-\good, the
case of the proposition that has already been
proved gives
\Equation\label{for large n}
\rank\pi_1(M)\le1+\bigg\lfloor\frac{\vol M}{b(\epsilon_n/2)}\bigg\rfloor\cdot
\max\bigg(0,
\frac{15}{32}\bigg\lfloor\frac{B(R)-b(\epsilon_n/2)}c\bigg\rfloor-1+\frac{\rho}{16}\bigg)
\EndEquation
for sufficiently large $n$. 
Taking limits in (\ref{for large n}) as
$n\to\infty$, and using the continuity of the function $b$ 
(see \cite[Subsection 4.1]{ratioI})
and the
upper semicontinuity of the greatest integer function $x\to\lfloor
x\rfloor$, we obtain the conclusion of the proposition. 
\EndProof

\Definition\label{recall}
As was mentioned in the introduction, a group $\Pi$ is
said to be {\it $k$-free} for a given positive integer $k$ if each 
subgroup of $\Pi$ having rank at most $k$ is free; and a group $\Pi$ is
said to be {\it $k$-semifree} for a given  $k$ if each 
subgroup of $\Pi$ having rank at most $k$ is a free product of free
abelian groups. 
\EndDefinition

\Corollary\label{almost main}
Let $M$ be a \willbefinitevolume\  orientable hyperbolic $3$-manifold
such that $\pi_1(M)$ is $2$-semifree. 
Suppose that $\rho$ is a non-negative integer such that, for every
point $v\in M$, the subgroup of $\pi_1(M,v)$ generated by all elements
represented by closed paths of length less than $\log9=2\log3$ based at
$v$ has rank at most $\rho$. Then we have
$$\rank\pi_1(M)\le1+\bigg\lfloor\frac{\vol M}{b((\log3)/2)}\bigg\rfloor\cdot\bigg(146.1875+\frac{\rho}{16}\bigg).$$
\EndCorollary

\Proof
We apply Proposition \ref{betwixt and between}, taking
$\epsilon=\log3$,  $R=2\log3+0.15$, and $c=0.496$. According to \cite[Corollary 4.2]{acs-surgery}, the hypothesis that
$\pi_1(M)$ is $2$-\willbesemifree\  implies that $\epsilon:=\log3$ is a Margulis
number for $M$. According to \cite[Lemma 5.1]{ratioI}, with the values
of $\epsilon$, $R$ and $c$ that we have specified, we have $\phi(R-D,\epsilon/2,D)>c$ for every $D$
in the interval
$[R/2-\epsilon/4,\epsilon]
=[3(\log3)/4+.075,\log3]$. Hence the conclusion of
Proposition \ref{betwixt and between} holds with these values of
$\epsilon$, $R$ and $c$. The quantity
$(B(R)-b(\epsilon/2))/c$, which appears in the conclusion of
Proposition \ref{betwixt and between}, is equal to
$314.62\ldots$. Since the quantities $(15/32)\times314-1$ and $\rho/16$ are non-negative, the conclusion of
Proposition \ref{betwixt and between} becomes
$$\rank\pi_1(M)\le1+\bigg\lfloor\frac{\vol M}{b((\log3)/2)}\bigg\rfloor\cdot\bigg(\frac{15}{32}\times314-1+\frac{\rho}{16}\bigg),$$
which gives the conclusion of the corollary.
\EndProof

\section{Proofs of the main results}

The main results of this paper are Theorems \ref{mainer},
\ref{semimainer}, \ref{closed homology} and \ref{maybe cusps}
below. The statements of these results were previewed in the introduction.

\Definition\label{indy def}
We shall say that elements $x_1,\ldots,x_m$ of a group $\redGamma $ are {\it independent} if the subgroup $\langle x_1,\ldots,x_m\rangle$ of $\redGamma $ is free on its generators $x_1,\ldots,x_m$.
\EndDefinition

The following result
will play an important role in this section.
As was mentioned in the introduction, it is a version of the ``
$\log(2k-1)$ Theorem.'' (The case $k=2$ of the theorem has already
been used indirectly in the paper, via the application of
\cite[Corollary 4.2]{acs-surgery} in the proof of Corollary
\ref{almost main}.)

\Theorem\label{another log(2k-1)}
Let $p$ be a point of an orientable hyperbolic $3$-manifold $M$,  let $k\ge2$ be
an integer, and let $x_1,\ldots,x_k$ be independent elements of
$\pi_1(M,p)$. For $i=1,\ldots,k$, let $d_i$ denote the minimum length
of a closed path based at $p$ and representing the element $x_i$. Then we have
$$
\sum_{i=1}^k\frac1{1+e^{d_i}}\le\frac12.
$$
In particular we have $d_i\ge\log(2k-1)$ for some
$i\in\{1,\ldots,k\}$.
\EndTheorem

\Proof
Write $M=\HH^3/\Gamma$, where $\Gamma\le\isomplus(\HH^3)$ is discrete
and torsion-free, and let $q:\HH^3\to M$ denote the quotient
map. Choose a point $P\in q^{-1}(p)$. By covering space theory, $P$
defines an isomorphism $j:\pi_1(M,p)\to\Gamma$; we set
$\xi_i=j(x_i)$ for $i=1,\ldots,k$. Then
$\xi_1,\ldots,\xi_k$ are independent elements of $\Gamma$, and
therefore freely generate a discrete subset of $\isomplus(\HH^3)$.
For $i=1,\ldots,k$ we have $d_i=\dist(P,\xi_i\cdot P)$. The assertions
now follow from \cite[Theorem 4.1]{acs-surgery}, with $P$ playing the
role of $z$ in the latter result.
\EndProof

\Proposition\label{independent}
Let $F$ be a free group of finite rank $r$, and let $S$ be a finite
generating set for $F$. Then $S$ contains $r$ independent elements.
\EndProposition

\Proof
Let $p$ denote the canonical
homomorphism from $F$ to $H_1(F,\QQ) $. Then $p(S)$ spans the
$\QQ$-vector space $H_1(F,\QQ)$, which has dimension $r$; hence $p(S)$ contains a basis of $H_1(F,\QQ)$,
which may  be written as $\{p(x_1),\ldots,p(x_r)\}$, where
$x_1,\ldots,x_r$ are elements of $S$. Set $E=\langle
x_1,\ldots,x_r\rangle\le F$, let $q$ denote the canonical
homomorphism from $E$ to $H_1(E,\QQ) $, and let $h: H_1(E,\QQ)\to
H_1(F,\QQ) $ denote the linear map induced by the inclusion $E\to
F$. Then $h\circ q=p|E$, and since $\{p(x_1),\ldots,p(x_r)\}$ is  a
basis for $H_1(F,\QQ)$, the map $h$ must be surjective. Hence $\dim
H_1(E;\QQ)\ge\dim H_1(F;\QQ)$, i.e. $\rank E\ge\rank F=r$. Since
$x_1,\ldots,x_r$ generate $E$, the rank of $E$ is exactly $r$;
furthermore, since the $r$ elements $x_1,\ldots,x_r$ generate the
rank-$r$ free group $E$, they must form a basis of $E$ 
(see  \cite[vol. 2, p. 59]{kurosh}).
Thus $x_1,\ldots,x_r$ are independent elements of $S$.
\EndProof

\Proposition\label{k-free consequence}
Let $k$ be a positive integer, and let $\redGamma $ be a $k$-free
group (see Definition \ref{recall}). Then for any finite set $S\subset \redGamma $,
either (a) $\rank \langle S\rangle<k$, or (b) $S$ contains $k$
independent elements.
\EndProposition

\Proof
Suppose that (a) does not hold, i.e. that $\rank \langle S\rangle\ge
k$. Among all subsets of $S$ that generate subgroups of $\redGamma $ having
rank at least $k$, choose one, $S_0$, which is minimal with respect to
inclusion. Then the subgroup $F:=\langle S_0\rangle$ of $\redGamma $ has
rank $k$. Since $\redGamma $ is $k$-free, $F$ is free. Applying
Proposition \ref{independent}, with $S_0$ playing the role of $S$, we deduce that $S_0$ contains $k$ independent
elements. In particular (b) holds.
\EndProof

\Proposition\label{semifree consequence}
Let $k$ be a positive integer, and let $\redGamma $ be a $(2k-1)$-semifree
group (see Definition \ref{recall}). Suppose that every free abelian subgroup of $\redGamma $ has rank at
most $2$. Then for any finite set $S\subset \redGamma $,
either (a) $\rank \langle S\rangle\le2k-2$, or (b) $S$ contains $k$
independent elements.
\EndProposition

\Proof
Suppose that (a) does not hold, i.e. that $\rank \langle S\rangle\ge
2k-1$. Among all subsets of $S$ that generate subgroups of $\redGamma $ having
rank at least $2k-1$, choose one, $S_0$, which is minimal with respect to
inclusion. Then the subgroup $E:=\rank \langle S\rangle$ of $\redGamma $ has
rank $2k-1$. Since $\redGamma $ is $(2k-1)$-semifree, $E$ is a free product of non-trivial
free abelian groups, say $A_1\star\cdots\star A_m$. If we set
$r_i=\rank A_i$ for $i=1,\ldots,m$, the hypothesis implies that
$r_i\le2$ for each $i$. We have $2k-1=\rank E=\sum_{i=1}^mr_i\le2m$,
and hence $m\ge k$. Since the $A_i$ are non-trivial free abelian
groups, $E$ has a quotient $F$ which is free of rank $k$. Let $q:E\to
F$ denote the quotient homomorphism. Then $q(S)$ generates $F$, and we may apply
Proposition \ref{independent}, with $q(S)$ playing the role of $S$ in
that proposition, to deduce that $q(S)$ contains $k$ independent
elements, say $y_1,\ldots,y_k$. For $i=1,\ldots,k$, we may write
$y_i=q(x_i)$ for some $x_i\in S$. Since $y_1,\ldots,y_k$ are
independent, $x_1,\ldots,x_k$ are also independent, and hence  (b) holds.
\EndProof

\Theorem\label{mainer}
Let $M$ be a \willbefinitevolume\  orientable hyperbolic $3$-manifold
such that $\pi_1(M)$ is $5$-free. 
Then we have
$$\rank\pi_1(M)<1+
157.497\cdot
\vol(M)
.$$
\EndTheorem

\Proof
Since $\pi_1(M)$ is $5$-free, it is in particular $2$-free and
therefore $2$-semifree. 

For each point $v\in M$, let $S_v$ denote the
set of  all elements of $\pi_1(M,v)$
represented by closed paths of length less than $\log9$ based at
$v$. For each  $v\in M$, 
since $\pi_1(M,v)$ is
$5$-free, it follows from Proposition \ref{k-free consequence} that either (a) $\rank \langle S_v\rangle<5$, or (b) $S_v$ contains five
independent elements. But Alternative (b) would contradict the case
$k=5$ of Theorem \ref{another log(2k-1)}.
Hence (a) holds for every $v\in M$. We may therefore apply Corollary
\ref{almost main}, taking $\rho=4$, to deduce that 
$$
\begin{aligned}
\rank\pi_1(M)&\le1+\bigg\lfloor\frac{\vol M}{b((\log3)/2)}\bigg\rfloor\cdot\bigg(146.1875+\frac{1}{4}\bigg)\\
&\le1+\frac{146.1875+1/4}{b((\log3)/2)}\cdot\vol(M)\\
&<1+157.497\cdot
\vol(M)
.
\end{aligned}
$$
\EndProof

\Theorem\label{semimainer}
Let $M$ be a \willbefinitevolume\  orientable hyperbolic $3$-manifold
such that $\pi_1(M)$ is $9$-semifree. 
Then we have
$$\rank\pi_1(M)<1+
157.766
\cdot
\vol(M)
.$$
\EndTheorem

\Proof
Since $\pi_1(M)$ is $9$-semifree, it is in particular $2$-semifree. 

For each point $v\in M$, let $S_v$ denote the
set of  all elements of $\pi_1(M,v)$
represented by closed paths of length less than $\log9$ based at
$v$. For each $v\in M$,
since $\pi_1(M,v)$ is
$9$-semifree, it follows from Proposition \ref{semifree consequence} that either (a) $\rank \langle S_v\rangle<9$, or (b) $S_v$ contains five
independent elements. But Alternative (b) would contradict 
the case
$k=5$ of Theorem \ref{another log(2k-1)}.
Hence (a) holds for every $v\in M$. We may therefore apply Corollary
\ref{almost main}, taking $\rho=8$, to deduce that 
$$
\begin{aligned}
\rank\pi_1(M)&\le1+\bigg\lfloor\frac{\vol M}{b((\log3)/2)}\bigg\rfloor\cdot\bigg(146.1875+\frac{1}{2}\bigg)\\
&\le1+\frac{146.1875+1/2}{b((\log3)/2)}\cdot\vol(M)\\
&<1+157.766\cdot
\vol(M)
.
\end{aligned}
$$
\EndProof


\Theorem\label{closed homology}
Let $M$ be a closed, orientable hyperbolic
$3$-manifold, and let $p$ be  a prime. Suppose that either 
\begin{enumerate}[(a)]
\item $\pi_1(M)$ has no subgroup isomorphic to the fundamental group
  of a closed, orientable  surface of genus $2$, $3$ or $4$; or
\item $p=2$, and $M$ contains no (embedded, two-sided) incompressible
  surface of genus $2$, $3$ or $4$. 
\end{enumerate}
Then 
$$\dim H_1(M;\FF_p)< 157.763\cdot\vol (M).$$
\EndTheorem

\Proof
First consider the case in which 
$\dim H_1(M;\FF_p)\le10$. 
According to 
\cite[Theorem 1.3]{milley},
we have $\vol M>0.94$ for any \willbefinitevolume\ orientable hyperbolic
$3$-manifold $M$. Hence in this case, for any prime $p$, we have $\dim
H_1(M;\FF_p)<11\cdot \vol (M)$, which is  stronger than the conclusion of the
theorem. For the rest of the proof, we will assume
$\dim H_1(M;\FF_p)\ge11$.

According to \cite[Theorem 1.1]{singular-two}, if $g\ge2$ is an integer, and if $M$
is a closed, orientable, hyperbolic $3$-manifold such that $\dim
H_1(M;\FF_2)
\ge
\max(3g - 1, 6)$, and if $\pi_1(M)$ has a
subgroup isomorphic to the fundamental group of a closed, orientable
surface  of genus  $g$, then M contains a closed, incompressible
surface of genus at most $g$. Hyperbolicity implies that the genus of
such a surface is at least $2$. Hence, under the assumption that 
$\dim
H_1(M;\FF_2)\ge11$, 
the existence of a
subgroup isomorphic to the fundamental group of a closed, orientable
surface  of genus  $2$, $3$ or $4$ implies the existence of a closed, incompressible
surface of genus $2$, $3$ or $4$. This means that, under the
assumption that $\dim H_1(M;\FF_p)\ge11$, Alternative (b) of the
hypothesis implies Alternative (a). We may therefore assume for the
rest of the proof that Alternative (a) holds.

We now apply \cite[Corollary 7.4]{accs}, which includes the assertion that if $k$ is
a positive integer, and if the
fundamental group of an orientable $3$-manifold $M$ has no
subgroup isomorphic to the fundamental group of any closed orientable
surface whose genus lies in the set  $\{2,3,\ldots, k - 1\}$, and if
$\dim H_1(M;\FF_p)\ge k + 2$ for some prime $p$, then
$\pi_1(M)$ is $k$-semifree; furthermore, according to 
\cite[Remark 7.5]
{accs}, 
in the case where $M$ is closed, if $\pi_1(M)$ is
$k$-semifree then it is $k$-free. Since Alternative (a) of the
hypothesis holds, and since $\dim H_1(M;\FF_p)\ge 11>7$, it now follows
that $\pi_1(M)$ is $5$-free. 
If we set $V=\vol M$, then by Theorem \ref{mainer} we have
\Equation\label{oliounidiselove}
\rank\pi_1(M)<1+
157.497\cdot V=\bigg(\frac1V+157.497\bigg)V.
\EndEquation
But since the closed, orientable hyperbolic $3$-manifold $M$ has a
$5$-free fundamental group, it follows from \cite[Proposition
12.13]{kfree-volume} that $V>3.77$. Hence
$$\frac1V+157.497<
\frac1{3.77}
+157.497=157.7622\ldots,$$
which with (\ref{oliounidiselove}) implies the conclusion of the theorem.

\EndProof

For the case where the manifold $M$ is not assumed to be closed, we
have the following result.

\Theorem\label{maybe cusps}
Let $M$ be a finite-volume, orientable hyperbolic
$3$-manifold, and let $p$ be  a prime. Suppose that for
$g=2,3,\ldots,8$, the group $\pi_1(M)$ has no subgroup isomorphic to the fundamental group
  of a closed, orientable  surface of genus $g$.
Then 
$$\dim H_1(M;\FF_p)< 
158.12\cdot\vol( M).$$
\EndTheorem

\Proof
If $M$ is closed, this follows from Theorem \ref{closed homology}. We
will therefore assume that $M$ is non-compact.

Consider the case in which 
$\dim H_1(M;\FF_p)\le10$. 
Since
\cite[Theorem 1.3]{milley} gives $\vol M>0.94$, for any prime $p$, we then have $\dim
H_1(M;\FF_p)<11\cdot V$, which is  stronger than the conclusion of the
theorem. 

For the rest of the proof, we will assume
$\dim H_1(M;\FF_p)\ge11$. Using that $\pi_1(M)$ has no
subgroup isomorphic to the fundamental group of any closed orientable
surface whose genus lies in the set  $\{2,3,\ldots, 8\}$, and 
applying the case $k=9$ of \cite[C
orollary 7.4]{accs}, we deduce that
$\pi_1(M)$ is $9$-semifree. 
If we set $V=\vol M$, then by Theorem \ref{semimainer} we have
\Equation\label{oioudoudouiou}
\rank\pi_1(M)<1+
157.766
=\bigg(\frac1V+157.766\bigg)V.
\EndEquation
But since the  orientable hyperbolic $3$-manifold $M$ is non-compact
and $\dim H_1(M;\FF_p)\ge11>3$, it follows from \cite[Lemma
5.3]{ratioI} that $V>2.848$. Hence
$$\frac1V+157.766<\frac1{2.848}+157.766=158.117\ldots,$$
which with (\ref{oioudoudouiou}) implies the conclusion of the theorem.
\EndProof

\bibliographystyle{plain}

\end{document}